\let\LaTeXcline\cline
\documentclass[sn-mathphys,Numbered]{sn-jnl}
\let\cline\LaTeXcline

\usepackage{graphicx}
\usepackage{algorithm, algorithmic}
\usepackage{float}

\usepackage{amssymb,latexsym,amsmath,enumerate,verbatim,amsfonts,mathtools}
\usepackage{amsthm}%
\usepackage{mathrsfs}%
\usepackage{manyfoot}
\usepackage{xcolor}%
\usepackage{booktabs}%
\usepackage{listings}%
\usepackage{mathptmx}

\usepackage[title]{appendix}%

\usepackage{bm,bbm}
\usepackage{caption}
\usepackage{subcaption}
\usepackage{physics}
\usepackage{multirow}
\usepackage{tikz}
\usetikzlibrary{positioning}
\usepackage{makecell}
\usepackage{hyperref}
\newtheorem{theo}{Theorem}
\newtheorem{lem}[theo]{Lemma}
\newtheorem{remark}[theo]{Remark}
\newtheorem{asmp}[theo]{Assumption}
\newtheorem{prop}[theo]{Proposition}
\newtheorem{alg}[theo]{Algorithm}

\DeclareMathOperator*{\argmin}{\arg\min}
\newcommand{\E}[1]{\mathbb{E}\left[#1\right]}
\newcommand{\pr}{\rm proj}

\usepackage{ifthen}
\newcommand{\COMM}[2]{{
\ifthenelse{\equal{#1}{AT}}{\color{red}}{
\ifthenelse{\equal{#1}{RK}}{\color{blue}}}
[#1: #2]
}}

\begin{document}

\title[Convergence error analysis of RGLD]{Convergence error analysis of reflected gradient Langevin dynamics for non-convex constrained optimization
}
\author[1]{\fnm{Kanji} \sur{Sato}}\email{descartes329@gmail.com}

\author[1,3]{\fnm{Akiko} \sur{Takeda}}\email{takeda@mist.i.u-tokyo.ac.jp}

\author*[2]{\fnm{Reiichiro} \sur{Kawai}}\email{raykawai@g.ecc.u-tokyo.ac.jp}

\author[1,3]{\fnm{Taiji} \sur{Suzuki}}\email{taiji@mist.i.u-tokyo.ac.jp}

\affil[1]{\orgdiv{Graduate School of Information Science and Technology}, \orgname{The University of Tokyo}, \orgaddress{\street{Bunkyo-ku}, \city{Tokyo}, \postcode{113-8656}, \country{Japan}}}

\affil*[2]{\orgdiv{Graduate School of Arts and Sciences / Mathematics and Informatics Center}, \orgname{The University of Tokyo}, \orgaddress{\street{Meguro-ku}, \city{Tokyo}, \postcode{153-8902},  \country{Japan}}}

\affil[3]{\orgdiv{Center for Advanced Intelligence Project}, \orgname{RIKEN}, \orgaddress{\street{Chuo-ku}, \city{Tokyo}, \postcode{103-0027}, \country{Japan}}}




\abstract{
    Gradient Langevin dynamics and a variety of its variants have attracted increasing attention owing to their convergence towards the global optimal solution, initially in the unconstrained convex framework while recently even in convex constrained non-convex problems. In the present work, we extend those frameworks to non-convex problems on a non-convex feasible region with a global optimization algorithm built upon reflected gradient Langevin dynamics and derive its convergence rates. By effectively making use of its reflection at the boundary in combination with the probabilistic representation for the Poisson equation with the Neumann boundary condition, we present promising convergence rates, particularly faster than the existing one for convex constrained non-convex problems.
}



\maketitle

\section{Introduction}

In the present work, we are interested in the following (generally) non-convex problems over a closed set $\overline K\subset\mathbb{R}^d$, defined as 
\begin{equation}\label{eq:problem}
\min_{x\in\mathbb{R}^d} f(x)  \quad \mathrm{s.t.}\quad x\in \overline K,
\end{equation}
where $f$ is a non-convex function mapping from ${\mathbb{R}}^d$ to $\mathbb{R}$ and $K$ is a (possibly non-convex) bounded open
region with a smooth boundary (whose assumptions are detailed later in   
Section \ref{sec:assumptions}).

The most common approach to such constrained optimization problems with boundaries has long been based on the projected gradient descent (PG), the procedure of which can be described as 
\[
    X_{k+1} = \mathcal{P}_K\left(X_k-\eta\grad f(X_k)\right),
\]
where $\mathcal{P}_K$ represents a projection operator defined later in \eqref{def:projection}.
In general, however, PG may well get stuck at a non-global stationary point.
In order to lead the algorithm out of such a point, it has been found effective to add stochastic perturbations, boosting up
gradient Langevin dynamics (GLD) to one of the most promising algorithms \cite{Raginsky,Welling}.

A line of research on Langevin dynamics, originating from stochastic diffusion, has been focused on global optimization of non-convex problems 
\cite{Chiang,Gelfand,geman1986diffusions}.
Constrained non-convex problems have not been tackled until recently because of theoretical difficulties of the underlying stochastic differential equations (SDE), whose solutions cannot be shown to exist in a trivial manner \cite{lions1984stochastic,Tanaka}, whereas recently in \cite{Lamperski}, projected gradient Langevin dynamics (PGLD) \cite{Bubeck} is introduced to address convex constrained non-convex problems, 
yet with a wide gap remaining in the
convergence rate for constrained and unconstrained problems.
We summarize the relevant convergence rates in Table \ref{table:related-works} (as well as relevant studies on GLD later in Appendix \ref{sec:related-work}).

\begin{table}[ht]
    \centering
    \begin{tabular}{|c|c|c|c|} \hline
        \multicolumn{2}{|c|}{}& \multicolumn{2}{c|}{objective function}\\ \cline{3-4}
        \multicolumn{2}{|c|}{}& convex & non-convex\\ \hline 
        \multicolumn{2}{|c|}{unconstrained} & $\widetilde{\mathcal{O}}(\epsilon^{-2})$ \cite{dalalyan2017further} & $\widetilde{\mathcal{O}}(d^4\epsilon^{-2})$ \cite{zou2021faster}\\ \hline
        \multirow{2}{*}{\makecell{constrained}} 
        & convex & $\widetilde{\mathcal{O}}(d\epsilon^{-2})$ \cite{hsieh2018mirrored} & $\widetilde{\mathcal{O}}(d^4\lambda_*^{-1}\epsilon^{-4})$ \cite{Lamperski}\\ \cline{2-4}
                                     & non-convex & \multicolumn{2}{c|}{\makecell{$\widetilde{\mathcal{O}}(d^3\lambda_*^{-3}\epsilon^{-3})$} (The present work)}\\ \hline
\end{tabular}
\caption{Iteration complexity of GLD for the convergence towards the stationary distribution, where $\widetilde{\mathcal{O}}$ denotes the order ignoring polylogarithmic factors, while $\lambda_*$ is the spectral gap defined in Section \ref{subsection main result}.
Here, the convergence is defined in different yet comparable ways.
The present work derives an upper bound for the convergence of $\mathbb{E}[N^{-1}\sum_{k=0}^{N-1} f(X_k)]-\min_{x \in \overline K} f(x)$, that ensures the inequality $\E{\min_k f(X_k)}- \min_{x \in \overline K} f(x)\leq \epsilon$, while 
this also yields a bound for $\E{f(X_{k'})}-\min_{x \in \overline K} f(x)$ with $X_{k'}$ a randomly sampled state from the discretized trajectory.
}
    \label{table:related-works}
\end{table}

The aim of the present work is to develop and investigate an algorithm with a smaller gap in the convergence rate for constrained and unconstrained problems.
To this end, we employ the following {\it reflected gradient Langevin dynamics (RGLD)}: 
\[
         X_{k+1}= \mathcal{R}_K\left(X_k - \eta\grad f(X_k) + \sqrt{\frac{2\eta}{\beta}}\xi_{k+1}\right),
\]
as an optimization algorithm for minimizing the non-convex objective function $f$ on a possibly non-convex feasible region $\overline K$.
Here, we denote by $\mathcal{R}_K$ the reflection operator defined later in \eqref{def:projection} and by $\{\xi_k\}_{k\in\mathbb{N}}$ a sequence of independent and identically distributed (iid) Rademacher random vectors in $\mathbb{R}^d$ (that is, their $d$ components are iid Rademacher random variables), while $\eta$ is a positive step size, and $\beta$ denotes an inverse temperature parameter.
In short, the RGLD algorithm is thus gradient Langevin dynamics (or unadjusted Langevin algorithm) with {\it reflection steps} introduced so that no intermediate steps leave the feasible region.
Although the marginal law $\mathcal{L}(X(t))$ of a continuous-time limit $\{X(t):\,t\ge 0\}$ (of the discrete-time algorithm $\{X_k:\,k\in\mathbb{N}_0\}$ above) converges in time to a unique stationary distribution (say, $\pi$) \cite{Bubeck} (that is, $\mathcal{L}(X(t))\to \pi$ as $t\to +\infty$),
such discrete-time algorithms have not been extensively explored yet in the framework of non-convex constrained non-convex problems.  
As such, in the present work, we present a convergence error analysis
on the discrete-time RGLD algorithm above in addressing the optimization problem \eqref{eq:problem}.

\vspace{1em}
\noindent {\bf Our contributions}

\noindent In short, the present work aims at enriching the existing theoretical analysis for constrained optimization problems, especially here with (generally non-convex) smooth constraints, on the basis of the RGLD algorithm. 

While a geometric ergodicity is derived in \cite{Lamperski} based on a coupling argument for the 1-Wasserstein metric on a convexity assumption of the feasible region $\overline K$, 
its time-discretization error is $\mathcal{O}((\eta \log k)^{1/4})$, slower than the unconstrained counterpart $\mathcal{O}(\eta^{1/2})$. 
This property is largely because PGLD is employed there with a rather rough estimate for the discretization error in terms of the convergence in mean.
Despite PGLD has a convergence guarantee in convex-constrained problems (for instance, \cite{Bubeck}), it is known to induce a suboptimal rate.

To address this issue, we replace the projection step with the so-called ``reflection,''
which enables one to derive better convergence rates by canceling out certain terms of the error expansion. 
To avoid some theoretical difficulties incurred by the reflection step, we exploit the probabilistic representation for the Poisson equation with the Neumann boundary condition, corresponding to its continuous-time limit \cite{Leimkuhler}, which results in simpler proofs with the aid of the established results on the continuous-time dynamics.
In addition, we succeed to explicitly specify the parameter dependence of the error bound when deriving the convergence rate, unlike \cite[Theorem 4.2]{Leimkuhler}.
Since, for instance, a large inverse temperature $\beta$ needs to be set in the optimization procedure, it is of great importance to have control of those key parameters.

Our contribution is summarized as follows:
\begin{itemize}
\item
  This is the first study to address non-convex constrained problems with a certain smoothness assumption on the feasible region.
\item
  We develop a GLD algorithm based upon reflection steps and the Rademacher noise, unlike projection steps and the Gaussian noise in the existing GLD algorithms. 
\item
  The obtained convergence rate is, as summarized in Table \ref{table:related-works}, faster than that of \cite{Lamperski} owing to a sharper discretization error with a small gap still remaining in the convergence rate for constrained and unconstrained problems.
\end{itemize}
Finally, while employing existing techniques in deriving the convergence rate, we make several technical advances as well, in relation to the reflected diffusion (for instance, canceling-out of the second-order terms of the Taylor expansion), each of which is an essential yet quite unfamiliar piece in the relevant context.

\section{The algorithm: RGLD}

To describe the reflected gradient Langevin dynamics (RGLD) method in detail, we first develop the notation that is used throughout the paper. 
We let $\partial K$, 
$K^c$ and $\overline K$ denote the boundary, 
complement and closure of the set $K$, respectively. 
We denote by $K_{-r}$ the set of points outside $K$ whose distance from $\overline K$ is at most $r$, that is,
\begin{align*}
    K_{-r}:=\left\{x\in K^c:\,  \mathrm{dist}(x, K)  
 \leq r\right\},
\end{align*}
where $ \mathrm{dist}(z, S)  := \min_{y \in \overline S} \|y-z\|$ for the point $z$ and set $S$.
Let 
 $\mathcal{C}^n(K)$ denote the set of $n$-times continuously differentiable functions on the open set $K$ and let $\mathcal{C}^n(\overline{K})$ be the set of functions in $\mathcal{C}^n(K)$ all of whose derivatives of order up to $n$ have continuous extensions to the closure $\overline{K}$.

We let $\widetilde{\mathcal{O}}$ represent the order symbol ignoring the poly-log arithmetic factors.
Namely, for sequences $\{a_n\}_{n\in\mathbb{N}}$ and $\{b_n\}_{n\in\mathbb{N}}$, we mean by $a_n=\widetilde{\mathcal{O}}(b_n)$ that there exist positive constants $c_1$ and $c_2$ such that $|a_n|\leq c_1 |b_n|\log^{c_2} |b_n|$ for all $n$.
Similarly, for sequences of positive elements $\{a_n\}_{n\in\mathbb{N}}$ and $\{b_n\}_{n\in\mathbb{N}}$, we write $a_n \preceq b_n$ if there exists a positive constant $c$ such that $a_n \leq c b_n$ for all $n$, and write $a_n \sim b_n$ if $a_n \preceq b_n \land a_n \succeq b_n.$
Finally, we let $[n]$ denote the set $\{0, 1, \cdots, n-1\}$.
We say that $\xi$ is a Rademacher random vector in $\mathbb{R}^d$ if it consists of independent and identically distributed (iid) components, each of which only takes two values $\pm 1$ equally likely.
For a compact subset $\overline K$ of $\mathbb{R}^d$, 
we define the projection $\mathcal{P}_K$ and reflection $\mathcal{R}_K$, respectively, as follows: 
\begin{equation}
    \label{def:projection}
    \mathcal{P}_K(x) := \argmin_{y\in \overline K} \norm{y-x},\quad\mathcal{R}_K(x) := 2\mathcal{P}_K(x) - x.
\end{equation}


With the notation defined above, we are now ready to describe the RGLD method.
Throughout, we call the parameters $\eta$ and $\beta$ appearing in the algorithm, respectively, the step size and the inverse temperature parameter.

\begin{alg}[Reflected gradient Langevin dynamics (RGLD)]
\label{def:rgld}
Repeat the following procedure; for $k\in \mathbb{N}_0(:=\{0,1,\cdots\})$,
\[
  \begin{cases}
            X'_{k+1}=X_k - \eta\grad f(X_k) + \sqrt{2\eta/\beta}\xi_{k+1},\\
            X_{k+1}= \mathcal{R}_K\left(X'_{k+1}\right),  
  \end{cases}
\]
where $X_0 \in K$, $\eta > 0$, $\beta > 0$ and $\{\xi_k\}_{k\in\mathbb{N}}$ is a sequence of iid Rademacher random vectors in $\mathbb{R}^d$.
\end{alg}

If $\mathcal{P}_K$ is efficiently computable by some oracle and 
$\grad f$ is easily accessible by some oracle, the computation of each iteration is simple.
In this paper, we are interested in estimating how many iterations are required to obtain an $\epsilon$-optimal solution.
The reflection $\mathcal{R}_K$ in the algorithm above (or, equivalently, the projection $\mathcal{P}_K$ in light of the definitions \eqref{def:projection}) can be uniquely determined even in cases where the domain $K$ lacks convexity, provided that its boundary exhibits sufficient smoothness and the location $X_{k+1}'$ (prior to reflection) remains sufficiently close to the domain $K$ after overshooting it.
Let us stress the significance of employing the bounded Rademacher noise $\{\xi_k\}_{k\in\mathbb{N}}$ (over the unbounded Gaussian one) in maintaining the proximity of $X_{k+1}'$ to the domain when combined with sufficiently small $\eta$ and/or large $\beta$.

Conventionally, the projection (that is, $\mathcal{P}_K$, rather than $\mathcal{R}_K$) has been adopted to GLD, as in the projected GLD (PGLD) algorithm, as in 
\cite{Bubeck,Lamperski}.
Also, the Gaussian noise has naturally and thus often been employed for approximating the diffusion term appearing in the continuous-time limit of GLD \cite{bossy2004symmetrized}, whereas it can, rarely yet with a positive probability, produce so large outputs that the solution may jump too far away from the feasible region. 
The Rademacher noise that we employ in the algorithm above is bounded and thus keeps the trajectory (before projection or reflection) close to the feasible region.
It also turns out (Appendix \ref{sec:detailed-proof}) that the Rademacher noise here plays an essential role in achieving higher-order approximation error by canceling out the second-order term of the Taylor expansion in the course of the error analysis.


          {The proposed algorithm can be thought of as a time-discretization version of the continuous-time reflected gradient Langevin diffusion $\{X(t):\,t\ge 0\}$, given by
\begin{equation}
\label{eq:rsde}
    dX(t) = -\grad f(X(t))dt + \sqrt{2/\beta}d W(t) -v(X(t)) dL(t),
\end{equation}
where $\{W(t):\,t\ge 0\}$ is the standard $d$-dimensional Brownian motion, $v(\cdot)$ denotes the outward pointing normal defined on the boundary $\partial K$, and $\{L(t):\,t\ge 0\}$ is the boundary local time, which can be defined as $L(t):=\lim_{\epsilon\to 0}\epsilon^{-1}\int_0^t \mathbbm{1}(X(s)\in \partial \Omega_{\epsilon})ds,$ with $\partial \Omega_{\epsilon}:=\{x\in \overline{K}:\,{\rm dist}(x,\partial K)\le \epsilon\}$ (which is a neighborhood inside the domain, as opposed to the extension $K_{-r}$), where the limit here is well defined almost surely as well as in quadratic mean (see, for instance, \cite{mixedopt}).}
The first term in \eqref{eq:rsde} represents a drift descending along the direction of $-\grad f$, while the second one is the diffusion term that perturbs the gradient descent by Gaussian noise.
The last term $-v(X(t)) dL(t)$ governs reflection \cite{dupuis1999convex,lions1984stochastic,Tanaka} to keep the trajectory in the domain $K$
and obviously requires a non-trivial treatment when the dynamics is discretized.

It is known \cite{Bubeck} that under suitable conditions, the continuous-time RGLD \eqref{eq:rsde} converges weakly to a unique stationary distribution, which is the Gibbs distribution with the probability density function $\pi$, given by
\begin{equation}
    \label{eq:gibbs-measure}
    \pi(x):= 
    \frac{\exp(-\beta f(x))}{\int_K \exp(-\beta f(y))dy},
\end{equation}
defined on the domain $K$. 
The explicit expression of the Gibbs distribution here enables one to bound the discrepancy between the expectation of the objective with respect to the stationary distribution and the optimal value \cite{Raginsky}. 
As such, our strategy is to show a weak convergence of the discrete-time dynamics to the stationary distribution by bounding the discretization error with the aid of the probabilistic representation for the Poisson equation with the Neumann boundary condition.




\section{Convergence Analysis}

In this section, we present and discuss the main theorem of the present work, which is the convergence rate of the RGLD algorithm under suitable technical conditions. 


\subsection{Assumptions}
\label{sec:assumptions}

First, we prepare technical conditions on the domain $K$, the objective function $f$ and the optimization problem \eqref{eq:problem}.

\begin{asmp}
\label{assumption}
\leavevmode 
\begin{enumerate}
    \setlength{\parskip}{0cm}
    \setlength{\itemsep}{0cm}
    
    
    \item \label{asmp:interior-sphere} $K$ is 
      an open, bounded and connected subset of $\mathbb{R}^d$ that contains the origin;
    \item \label{asmp:boundary-smoothness}
    $\partial K\in \mathcal{C}^5$;
\item \label{asmp:func} $f\in \mathcal{C}^5(\overline{K})$;
    \item \label{asmp:solution-existence} The optimization problem \eqref{eq:problem} admits at least one optimal solution. 
\end{enumerate}
\end{asmp}


We impose Assumption \ref{assumption}-\ref{asmp:interior-sphere} to ensure the uniqueness of a stationary distribution of the reflected diffusion process \eqref{eq:rsde} (for instance, \cite{bencherif2009probabilistic,miranda1970}), which turns out to play a key role in the global optimization of RGLD.
It is worth noting that this base condition does not rule out non-convexity of the domain $K$.
Clearly, this condition also implies that the domain $K$ contains a Euclidean ball of a strictly positive radius (say, $r>0$) and also that $K$ is contained in a Euclidean ball with a positive radius (say, $R>0$) centered at the origin.
Note that these radii $r$ and $R$ appear in the statement of Lemma \ref{lem:near-optimality}, and that this assumption prevents the reflection operation from pushing the solution again out of the domain.


Assumption \ref{assumption}-\ref{asmp:boundary-smoothness} indicates, roughly speaking, that the boundary $\partial K$ can be thought of as the graph of such a smooth function.
This condition is imposed here so that an associated PDE (given shortly in \eqref{eq:neumann-problem}) admits a unique and smooth enough solution (for deriving Lemma \ref{lem:discretization-error}) and may seem strong, given that constrained problems are often described with multiple constraints and thus may have several non-differentiable extreme points.
For instance, Riemannian manifolds, such as Stiefel and Grassmann manifolds, meet the smoothness condition
\cite{patterson2013stochastic, wang2020fast}, whereas the $l_q$ norm ($q<1$), which is one of the most popular non-convex constraints for the sparse estimation, violates Assumption \ref{assumption}-\ref{asmp:boundary-smoothness} at points with zero elements.
Nonetheless, such constraints can be amended in such a way to have a smooth boundary with the aid of suitable smoothing operations.
After a smoothing operation, such as the Gaussian filter, has been applied as often done in practice, this smoothness condition is fulfilled even for non-smooth constraints appearing in the sparse estimation.

Finally, Assumption \ref{assumption}-\ref{asmp:func} implies that the objective function $f$ and the norm $\|\grad f\|$ are both uniformly bounded over the domain $\overline K$, since it is compact due to Assumption \ref{assumption}-\ref{asmp:interior-sphere}.
Moreover, the $M$-smoothness can easily be established, as below.

\begin{prop}
    \label{prop:m-smooth}
    Under Assumption \ref{assumption}, there exists a positive constant $M$ such that 
    $\norm{\grad f(x)-\grad f(y)} \leq M\norm{x-y}$ for all $x,y\in \overline K$.
\end{prop}

It is worth mentioning that, 
from a practical standpoint, the algorithm runs properly under Assumption \ref{assumption}-\ref{asmp:interior-sphere} alone, which is essential for the feasibility of the algorithm.

The reflection $\mathcal{R}_K$ in Algorithm \ref{def:rgld} can be uniquely determined if the step size $\eta$ is set sufficiently small, due to Assumption \ref{assumption}. 
For this, it suffices to show the asymptotics 
\begin{equation}
\left\|\mathcal{P}_K(X'_k)-X'_k\right\|=\mathcal{O}\left(\eta+\sqrt{d\eta/\beta}\right), 
\label{eq:proj_dist}
\end{equation}
as follows
\cite{bossy2004symmetrized,Leimkuhler}.
First, observe that $\norm{\mathcal{P}_K(X'_k)-X'_k}$ is bounded, based on the definition of the projection \eqref{def:projection} and $X_{k-1}\in K$, as 
\[
\norm{\mathcal{P}_K(X'_k)-X'_k} = \min_{x\in \overline K}\norm{x-X'_k} \leq \norm{X_{k-1}-X'_k},
\]
which can be further bounded as 
\[
\norm{X_{k-1}-X'_k} = \norm{- \eta\grad f(X_{k-1}) + \sqrt{2\eta/\beta}\xi_{k}}\leq \eta\norm{\grad f(X_{k-1})} + \sqrt{2\eta/\beta}\norm{\xi_k},
\]
due to the triangle inequality.
Since the norm $\|\grad f\|$ is uniformly bounded (by a suitable constant $G$) and $\norm{\xi_k}=\sqrt{d}$, it holds that
 \[
         \eta\norm{\grad f(X_{k-1})} + \sqrt{2\eta/\beta}\norm{\xi_k} \leq \eta G + \sqrt{2d\eta/\beta} = \mathcal{O}\left(\eta+\sqrt{d\eta/\beta}\right),
\]
which ensures that $\mathcal{R}_K(X'_k)$ does not overshoot the domain $K$, as long as the step size $\eta$ is sufficiently small.
In each experiment of Section \ref{chap:numerical-experiments}, we follow the usual practice of numerically verifying that every trajectory remains within the domain after each reflection, rather than adhering to the conservative requirement that the step size $\eta$ be smaller than the smallest radius of a finite open cover of the compact domain $\overline K$.


\subsection{Convergence rate of RGLD}
\label{subsection main result}

Here, we turn to the main result of the present work on the convergence rate of RGLD 
on the basis of the geometric ergodicity of the continuous-time limit of RGLD.
The rate of its weak convergence to the unique stationary distribution is governed by the spectral gap $\lambda_*$, defined as
the minimum non-zero eigenvalue of $-\mathcal{A}$, where $\mathcal{A}$ denotes the infinitesimal generator of the underlying dynamics, that is, for $h\in \mathcal{C}^2(K)$, 
\begin{equation}\label{infinitesimal generator}
\mathcal{A}h(x) := \beta^{-1}\Delta h(x) - \langle \grad f(x), \grad h(x) \rangle,\quad x\in K,
\end{equation}
where $\Delta$ denotes the Laplace operator.
For more detail on the spectral gap, we refer the reader to, for example, \cite{chen2012effects,peng2019asymptotic} as well as \cite[Chapter 4]{Bakry2014}.
We are now ready to state the main result.



\begin{theo}\label{thm:main}
Let Assumption \ref{assumption} hold, let $\eta < 1$ and let $\beta >1$.
Then, it holds that  
\[
        \E{\min_{k\in[N]}f(X_k)}- \min_{x \in \overline K} f(x) \preceq  
\frac{1}{\lambda_*N\eta}+\frac{\sqrt{\eta}\left(\beta\sqrt{\eta}+(\beta\eta+d)^{3/2}\right)}{\lambda_*}+\frac{d}{\beta}\log\beta.
\]
In particular, for every $0 < \epsilon \ll 1$, it holds that 
    \begin{equation}
        \E{\min_{k\in[N]}f(X_k)}- \min_{x \in \overline K} f(x) \preceq \epsilon+\frac{d\log\beta}{\beta},
\label{eq:error_bound}
    \end{equation}
    provided that
    \[
        \beta\succeq 1,\;
        \eta\preceq\min\left\{\frac{\lambda_*^2\epsilon^2}{d^3},\frac{\lambda_*\epsilon}{\beta}\right\},\;
        N\succeq \frac{1}{\lambda_*\epsilon\eta}\succeq \max\left\{\frac{d^3}{\lambda_*^3\epsilon^3}, \frac{\beta}{\lambda_*^2\epsilon^2}\right\}.
    \]
\end{theo}

Broadly speaking, the theorem claims that, by letting the step size $\eta$ be sufficiently small, an $\epsilon$-sampling error is achieved after $\widetilde{\mathcal{O}}(\max\{{d^3}{\lambda_*^{-3}\epsilon^{-3}}, {\beta}{\lambda_*^{-2}\epsilon^{-2}}\})$ iterations of RGLD.
Recall that the iteration complexity for PGLD obtained in \cite{Lamperski} is given as $N=\widetilde{\mathcal{O}}(d^4\lambda_*^{-1}\epsilon^{-4})$ (see Table~\ref{table:related-works}), which is not as sharp as our evaluation $N=\widetilde{\mathcal{O}}(d^3\lambda_*^{-3}\epsilon^{-3})$ 
when the inverse temperature parameter $\beta$ is sufficiently large in the order $\widetilde{\mathcal{O}}(d)$ (which is indeed the situation of interest). 

This improvement is mainly due to the reflection $\mathcal{R}_K$ employed in the present work, instead of the projection $\mathcal{P}_K$. 
With the reflection, the discretization error can be evaluated by approximating the time differentiation for the infinitesimal generator by a discrete-time differentiation with the aid of the Poisson equation with the Neumann boundary condition:
\begin{equation}
    \label{poisson neumann problem}
    \begin{cases}
         \mathcal{A}u(x) = f(x)-\mathbb{E}_\pi f, & x\in K,\\
         \langle \grad u(x), v(x)\rangle = 0, & x\in \partial K,
    \end{cases}
\end{equation}
where $\mathcal{A}$ is the infinitesimal generator \eqref{infinitesimal generator}.
In that approximation, the second-order derivatives are canceled out thanks to the symmetry between the states $X_{k}$ (after reflection) and $X'_{k}$ (before reflection), yielding a sharper bound than that of PGLD.
The discretization error between the outputs $f(X_k)$ (in discrete time) and $f(X(k\eta))$ (in continuous time) in PGLD is $\widetilde{\mathcal{O}}(\eta^{1/4})$ in \cite[Lemmas 5 and 6]{Lamperski}, not as sharp as our result. 




\begin{remark}{\rm
    The dependence of the spectral gap $\lambda_*$ on the inverse temperature parameter $\beta$ 
    has been studied in \cite{chen2012effects,peng2019asymptotic}, while it remains an open problem.
    One may be tempted to set the spectral gap $\lambda_*$ in proportion to the exponential of the inverse temperature parameter $\beta$, for the reason that the exponential dependence holds true for unconstrained GLD under similar assumptions \cite{bardet2018functional,gayrard2005metastability,Raginsky}.
    If that is really the case, finding an $\epsilon$-optimal solution suffers from the curse of dimensionality, where $\epsilon$-optimal solution $x$ is defined as the solution satisfying the bound $\E{f(x)}-\min_{x \in \overline K} f(x) \leq \epsilon$.
    It is necessary to set $\beta\succeq d\epsilon^{-1}\log\epsilon^{-1}$ to ensure $\mathbb{E}[\min_{k\in[N]}f(X_k)]-\min_{x \in \overline K} f(x) \leq \epsilon$, implying $N\succeq \exp (c_1d\epsilon^{-1})$.
    That is, the number of iterations $N$ depends on both the exponential of the problem dimension $d$ and that of (the inverse of) the error $\epsilon$.
\qed}\end{remark}


To derive Theorem \ref{thm:main},
we provide the following result.
 
\begin{lem}
    \label{lem:main}
    Let Assumption \ref{assumption} hold, let $\eta\preceq 1$ and let $\beta\succeq 1$.
    Then, there exists a positive constant $C$ such that for all $N\in \mathbb{N}$,
 \[
        \E{\frac{1}{N}\sum_{k\in [N]}f(X_k)}-\min_{x \in \overline K} f(x) 
        \le C\left( \frac{1}{\lambda_*N\eta}+\frac{\sqrt{\eta}\left(\beta\sqrt{\eta}+(\beta\eta+d)^{3/2}\right)}{\lambda_*}+\frac{d\log(C\beta)}{\beta}\right).
\]
\end{lem}

With this inequality in mind, one can easily find the parameter configuration to obtain a solution with $\epsilon$-sampling error by controlling the first two terms on the right-hand side. 
First, it is obvious from 
$\min_{k\in[N]} a_k \leq N^{-1}\sum_{k\in [N]}a_k$ 
that 
$$\mathbb{E}\left[\min_{k\in[N]}f(X_k)\right]-\min_{x \in \overline K} f(x)\leq \mathbb{E}\left[\frac{1}{N}\sum_{k\in [N]}f(X_k)\right]-\min_{x \in \overline K} f(x).$$
Thus, it suffices to show $\mathbb{E}[N^{-1}\sum_{k\in [N]}f(X_k)]-\min_{x \in \overline K} f(x)\preceq\epsilon + d\beta^{-1}\log\beta$ to derive 
\eqref{eq:error_bound}.
By Lemma \ref{lem:main}, it is enough to establish
\[
    \frac{1}{\lambda_*N\eta}+\frac{\sqrt{\eta}(\beta\sqrt{\eta}+(\beta\eta+d)^{3/2})}{\lambda_*}\preceq\epsilon.
\]
From the condition on $\eta$, we have
$\lambda_*^{-1}{\beta\eta}\preceq \epsilon$ and $\lambda_*^{-1}d^{3/2}\eta^{1/2}\preceq \epsilon$.
Moreover, we have $\lambda_*^{-1}\beta^{3/2}\eta^2\preceq \lambda_*\epsilon^2/\sqrt{\beta}$ and $\lambda_*\epsilon^2/\sqrt{\beta}\preceq\epsilon$, due to $\lambda_*,\epsilon,\beta^{-1/2}\preceq 1$.
Also, the condition on $N$ gives
$(\lambda_*N\eta)^{-1}\preceq \epsilon$.
Combining those yields Theorem \ref{thm:main}.

We hence turn to Lemma \ref{lem:main}.
Recall that $\{X_k:\,k\in\mathbb{N}_0\}$ denotes a time-discretized version of the continuous-time stochastic process $\{X(t):\,t\ge 0\}$ governed by the SDE \eqref{eq:rsde}. 
It is known (Lemma \ref{lem:estimates}) that the continuous-time limit $\{X(t):\,t\ge 0\}$ tends in distribution geometrically to the stationary distribution \eqref{eq:gibbs-measure}, 
that concentrates around an optimal solution $\argmin_{x\in \overline K} f(x)$. 
Thereby, it suffices to show that the discrete-time version $\{X_k:\,k\in\mathbb{N}_0\}$ well approximates the stationary distribution $\pi$ for sufficiently large $k$ and small $\eta$.
Along this line, we decompose the expected risk into two terms:
\[
\E{\frac{1}{N}\sum_{k\in [N]}f(X_k)}-\min_{x\in \overline K} f(x) = \left(\E{\frac{1}{N}\sum_{k\in [N]}f(X_k)}-\mathbb{E}_\pi f\right) + (\mathbb{E}_\pi f-\min_{x\in \overline K} f(x)),
\]
so that the following two lemmas provide an upper bound on the whole.
First, Lemma~\ref{lem:discretization-error} plays an important role for our purpose to show that the sequence generated by RGLD approximates sampling from the stationary distribution $\pi$ given a sufficiently small step size $\eta$ and a sufficiently large number of iterations $N$.

\begin{lem}
    \label{lem:discretization-error}
    Let Assumption \ref{assumption} hold, let $\eta<1$ and let $\beta<1$.
    Then, it holds that
    \begin{equation}
    \label{eq:lemdiscretization-error-bound}
        \E{\frac{1}{N}\sum_{k\in [N]}f(X_k)}-\mathbb{E}_{\pi}f\preceq \frac{1}{\lambda_*N\eta}+\frac{\sqrt{\eta}\left(\beta\sqrt{\eta}+(\beta\eta+d)^{3/2}\right)}{\lambda_*}.
    \end{equation}
\end{lem}

The next lemma, due to {\cite[Lemma 16]{Lamperski}, evaluates the quality of the expectation of the objective with respect to the stationary distribution $\pi$ as an approximation to the minimum value of the objective function.
Below, we let $R$ and $r$ denote strictly positive radii of outer and inner Euclidean balls of the feasible region $K$, respectively, owing to Assumption \ref{assumption}-\ref{asmp:interior-sphere}. 

\begin{lem}[\textcolor{black}{\cite[Lemma 16]{Lamperski}}]
    \label{lem:near-optimality}
    Under Assumption \ref{assumption}, there exists a positive constant $C$ such that
    \[
        \mathbb{E}_\pi f - \min_{x\in \overline K}f(x) \leq \frac{d}{\beta}\log(2R\max\left\{\frac{2}{r},\frac{C\beta(r+\sqrt{r^2+R^2})}{r\log{2}}\right\}).
    \]
\end{lem}

We note that \cite[Lemma 16]{Lamperski} is derived when the domain $K$ is convex.
In the present work, we have assumed that the domain $K$ is 
an open set (that is, $K$ is thick) and the optimal solution lies in $K$.
Since a ball of a sufficiently large radius can then be constructed around the optimal solution, this lemma can be employed in our context as well. 
Note that a similar formula can also be found in \cite{Raginsky} for unconstrained problems.





We close this section with a sketchy derivation of Lemma \ref{lem:discretization-error}.
Note that the existence and uniqueness of reflected SDEs \eqref{eq:rsde} have been studied in the name of the Skorokhod problem \cite{lions1984stochastic,Tanaka}.
Also, in a similar way to the unconstrained framework, the existence and uniqueness of a stationary distribution hold true for SDEs \eqref{eq:rsde} (see, for instance, \cite{Bubeck}).
For the sake of completeness, we provide a brief proof in Appendix \ref{sec:uniqueness}.

\begin{lem}
    \label{lem:unique-invariant-measure}
    Under Assumption \ref{assumption}, the stochastic process $\{X(t):\,t\ge 0\}$ admits a unique stationary distribution.
    Then, the stationary distribution is as defined in 
    \eqref{eq:gibbs-measure}.
\end{lem}

It thus remains to show the weak convergence of the time-discretized version and bound the error for a finite step $k$ and a positive step size $\eta$.
To this end, recall the Poisson equation with the Neumann boundary condition \eqref{poisson neumann problem}, which we restate here for the sake of convenience: 
\begin{equation}
    \label{eq:neumann-problem}
    \begin{cases}
         \mathcal{A}u(x) = f(x)-\mathbb{E}_\pi f, & x\in K,\\
         \langle \grad u(x), v(x)\rangle = 0, & x\in \partial K.
    \end{cases}
\end{equation}
As is known \cite{freidlin2016functional,miranda1970}, the so-called compatibility condition is required here:
\begin{align*}
    \int_K\left(f(x)-\mathbb{E}_\pi f\right)d\pi(x) = 0,
\end{align*}
which holds true by definition in our context.
{Thanks to Assumption \ref{assumption}, 
the Neumann problem above admits a unique solution (up to an additive constant), which is at least as smooth as the boundary condition or the objective function.
(We refer the reader to, for instance, \cite[Section 4.1]{Leimkuhler}.)}

In addition, in accordance with \cite[Theorem 4]{bencherif2009probabilistic}, the solution can be written as 
\begin{align*}
 u(x) = -\int_0^\infty \mathbb{E}\left[f(X^x(t))-\mathbb{E}_\pi f \right] dt,
\end{align*}
where $X^x(t)$ denotes the solution of the SDE \eqref{eq:rsde} at time $t$ with the initial state $x$. 
Then, by the Ito formula and the representation \eqref{eq:neumann-problem}, it holds that
$$\frac{\partial \mathbb{E}[u(X^x(t))]}{\partial t}\Big|_{t=0} = \mathcal{A} u(x) = f(x) - \mathbb{E}_\pi f.$$
By substituting $x = X_k$ here, we obtain an estimate of the difference $f(X_k) - \mathbb{E}_\pi f$ through the solution of the Poisson equation.
By further taking the average over $k\in [N]$, we get 
\[
\frac{1}{N}\sum_{k\in [N]} f(X_k) - \mathbb{E}_\pi f = \frac{1}{N}\sum_{k\in [N]}  \frac{\partial\mathbb{E}[u(X^{X_k}(t))]}{\partial t} \big|_{t=0}.
\]
Then, loosely speaking, we approximate each summand on the right-hand side as 
\[
\frac{\partial\mathbb{E}[u(X^{X_k}(t))]}{\partial t} \big|_{t=0} \approx \frac{1}{\eta}\left(\mathbb{E}\left[u(X_{k+1})\right] - \mathbb{E}\left[u(X_k)\right]\right),
\]
and thus arrive at the approximation:
\[
\frac{1}{N}\sum_{k\in [N]} f(X_k) - \mathbb{E}_\pi f \approx \frac{1}{N \eta} \left(\mathbb{E}\left[u(X_{N})\right] - \mathbb{E}\left[u(X_0)\right]\right).
\]

The right-hand side can be bounded by 
$\mathcal{O}(1/(\lambda_* N \eta))$ (due to Lemma \ref{lem:derivative-bounds} given in Appendix \ref{sec:detailed-proof}), corresponding to the first term of the upper bound \eqref{eq:lemdiscretization-error-bound}.
The second term of the upper bound \eqref{eq:lemdiscretization-error-bound} comes from the discretization error which can be evaluated as follows. 
In parallel to \cite[Theorem 4.2]{Leimkuhler} along with additional dependence of their bounds on the inverse temperature parameter $\beta$ and the problem dimension $d$, we perform the Taylor expansion to the increment $u(X_{k+1}) -u(X_k)$ up to the third order 
and bound the derivatives 
with the aid of the {\it Bismut-Elworthy-Li} formula  \cite{andres2011pathwise,cerrai1998some}.
In that evaluation, we show that 
the derivative of 
$\mathbb{E}\left[f(X^x(t))\right]-\mathbb{E}_\pi f$
with respect to the initial state $x$ decays exponentially as $t$ tends to infinity, 
along the line of \cite[Lemma 5.4]{Brehier}. 
Here, the reflection operation $\mathcal{R}_K$ plays a key role for a more precise evaluation of the discretization error in the sense of canceling the second-order derivative term. 
For a detailed derivation of Lemma \ref{lem:discretization-error}, we direct the reader to Appendix \ref{sec:detailed-proof}. 

\section{Numerical examples}
\label{chap:numerical-experiments}


In this section, we present numerical results to illustrate the effectiveness of the proposed RGLD algorithm on optimization problems with non-convex feasible regions, in comparison to PGLD (which has no theoretical convergence guarantee for such problems) and
with respect to the hyperparameters $\eta$ and $\beta$, as well as the problem dimension $d$.
In the experiments below, to ensure that the chosen values of $\eta$ and $\beta$ are valid, we have numerically verified that every trajectory remains within the domain after each reflection.
Throughout, we quantify the optimization error via the gap between objective function values at the global optimal solution and at the solution of the algorithms.



\subsection{Two-dimensional domain}
We start with a two-dimensional problem in order to present numerical results along with visualization and interpretation of the proposed algorithm in full detail.
The test function here is a Gaussian mixture density given by
\begin{equation}
    \label{gm}
    f(x) = -\sum_{k=1}^M w_k \exp(-\frac{1}{2}(x-m_k)^T(x-m_k)),
\end{equation}
where $\{w_k\}_{k\in \{1,\cdots,M\}}$ are randomly generated weights and $\{m_k\}_{k\in \{1,\cdots,M\}}$ are points with $M=5$ in the mesh grid $\{-2, -1, 0, +1, +2\}^2$.
The feasible region we define here lies between two spheres centered at the origin, with radii of $0.9$ and $4$, respectively,
which is clearly non-convex and satisfies Assumption \ref{assumption}.

\begin{figure}[ht]
  \centering
  \begin{minipage}[b]{0.45\linewidth}
    \centering
    \includegraphics[keepaspectratio, scale=0.4]{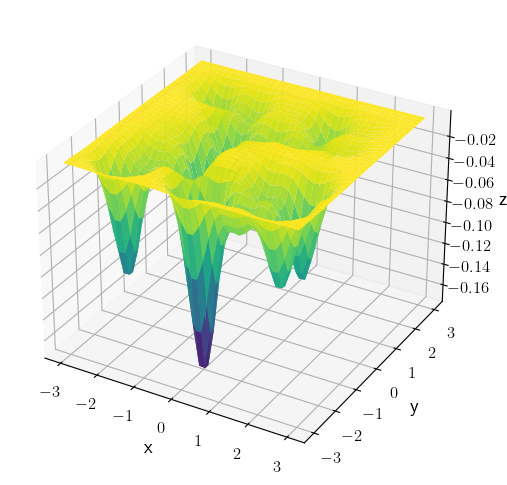}
    \subcaption{Objective function given in \eqref{gm}}
    \label{fig:gm-3dplot}
  \end{minipage}
  \begin{minipage}[b]{0.45\linewidth}
    \centering
    \includegraphics[keepaspectratio, scale=0.35]{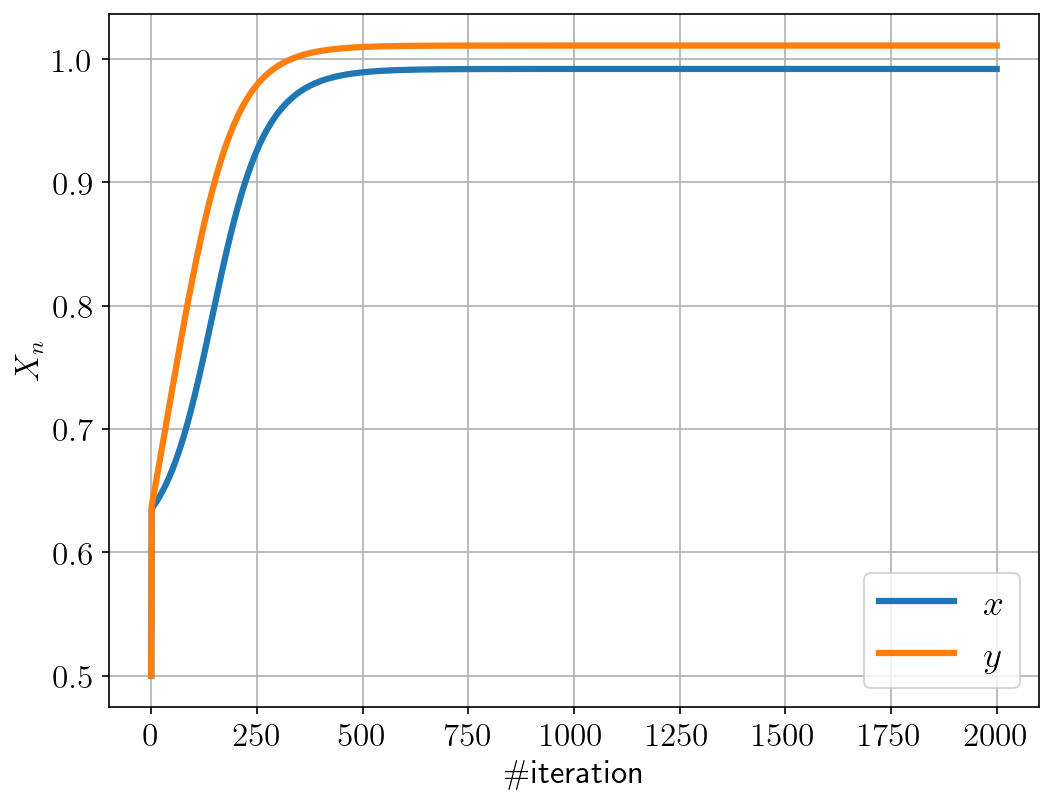}
    \subcaption{PG ($\eta=0.05$)}
     \label{fig:gm-error}
  \end{minipage}  \\
   \begin{minipage}[b]{0.45\linewidth}
    \centering
    \includegraphics[keepaspectratio, scale=0.35]{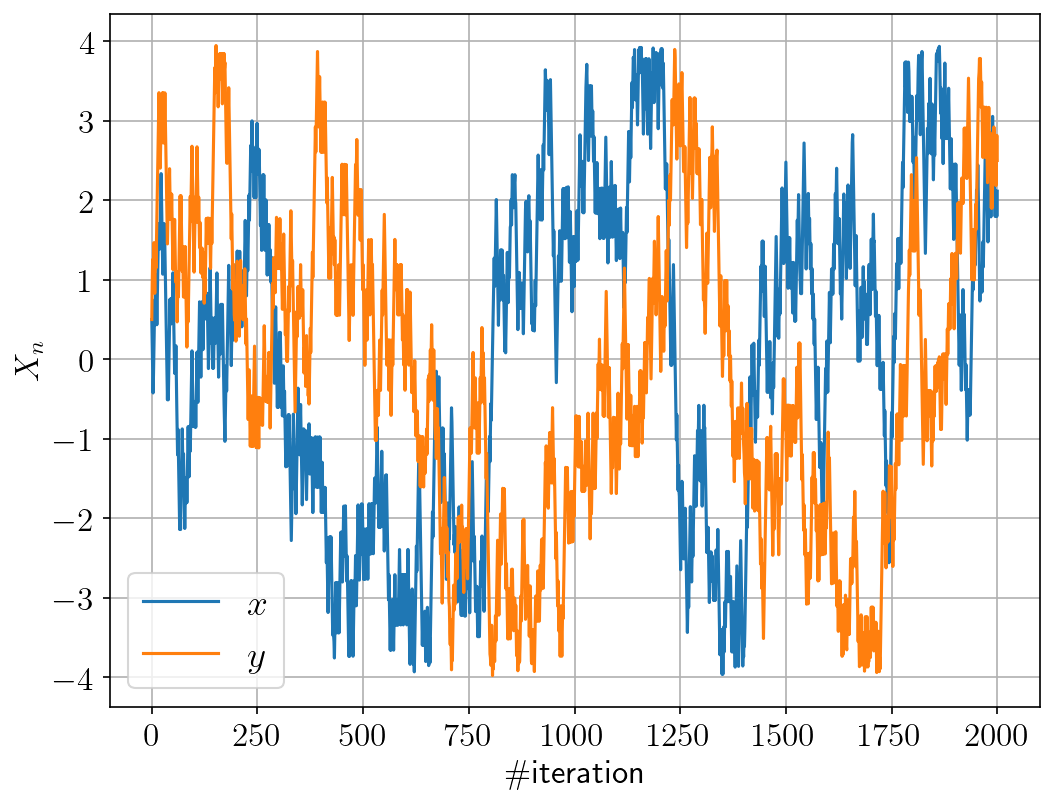}
    \subcaption{RGLD (\textcolor{black}{$\eta=0.05$ and} $\beta=1.0$)}
     \label{fig:gm-trajectory}
   \end{minipage}
    \begin{minipage}[b]{0.45\linewidth}
    \centering
    \includegraphics[keepaspectratio, scale=0.35]{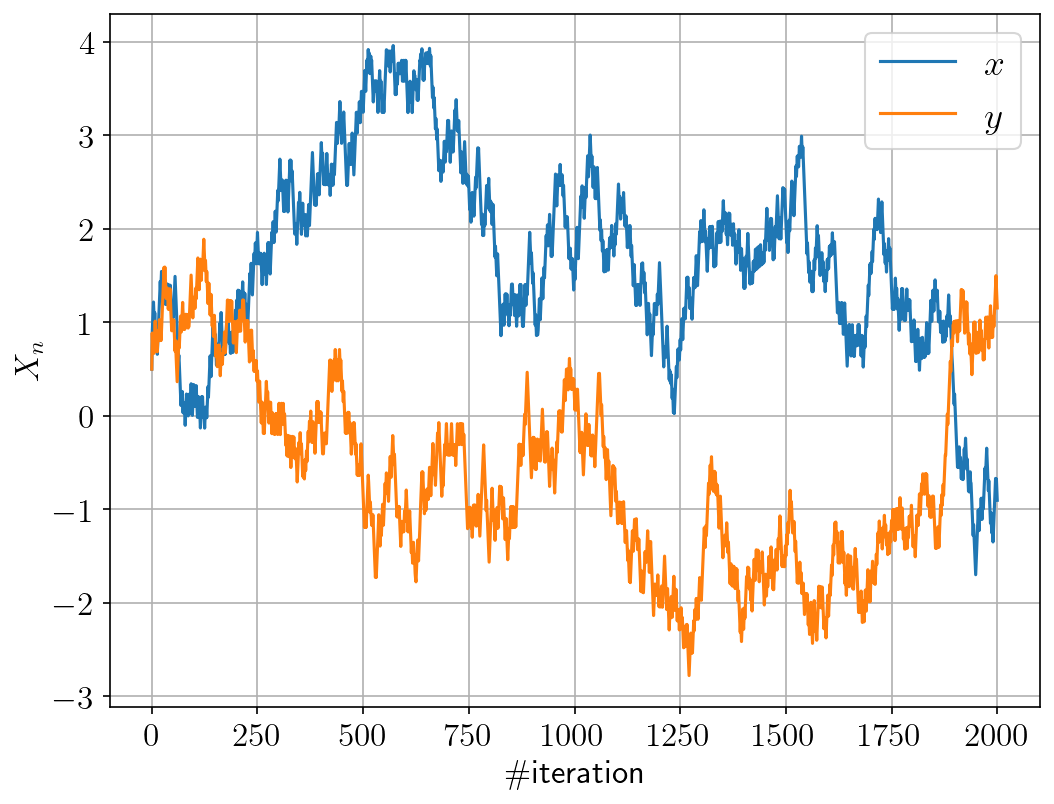}
    \subcaption{RGLD (\textcolor{black}{$\eta=0.05$ and} $\beta=8.0$)}
     \label{fig:gm-trajectory-betas}
  \end{minipage}
  \caption{(a) The objective function; (b)-(d) typical trajectories of the algorithms}
\end{figure}

In Figure \ref{fig:gm-3dplot}, we provide a 3D plot of the objective function with the global minimizer at $(0, -2)$ and 25 local minima.
The blue and orange lines in Figures \ref{fig:gm-error}, \ref{fig:gm-trajectory}, and \ref{fig:gm-trajectory-betas} correspond, respectively, to the first and second components of the solution $x$ in $\mathbb{R}^2$ found in each iteration of the respective algorithms.
To be more precise, Figure \ref{fig:gm-error} 
\textcolor{black}{depicts a typical trajectory}
of PG 
with $\eta=0.05$ 
and the initial state $X_0=(0.5, 0.5)$.
In this experiment, the PG algorithm is trapped at a local minimum and is not updated much anymore after $500$ steps or so. 
Next, in Figure \ref{fig:gm-trajectory}, we plot a typical trajectory of RGLD, which explores a wide range of the feasible region with rapid and frequent updates.
As illustrated in Figure \ref{fig:gm-trajectory-betas}, a larger $\beta$ stabilizes the trajectory, which may also help keep the iteration stuck in a mode of the objective function.
We note that a global optimal solution is found after around $1250$ steps (and then is left behind afterwards) if $\beta=8.0$ in this particular experiment.

The main result (Theorem \ref{thm:main}) gives an insight on how to set the hyperparameters $\eta$ and $\beta$ in an objective manner.
To examine the result, we plot the running minimum of the optimization error generated by the iterates of RGLD. 
Figures \ref{fig:gm-etas} and \ref{fig:gm-betas} present the convergence of RGLD for five different values of $\eta$ and $\beta$ (with the other one fixed), respectively.
A larger $\beta$ tends to restrict the solution to a narrower range of the region, resulting in slower convergences.
As for the parameter $\eta$, there seems to be a certain trade-off.
That is, $\eta$ needs to be set large enough so as to search a broad range of the region within a limited number of iterations.
In addition, the range has to be small enough to find narrow loss valleys as well.

It is also of great interest to compare the effects of the reflection and projection operations. 
Note that convergence rates of
PGLD are given in \cite{Bubeck,Lamperski}, 
but not for non-convex constrained problems, such as the one in the present experiment.
Figure \ref{fig:reflection-vs-projection} presents convergences of RGLD and PGLD.
Despite the reflection has been found theoretically preferable, this particular numerical experiment does not seem to demonstrate a significant difference in performance between those two operations. 

\begin{figure}[tb]
  \centering
   \begin{minipage}[b]{0.45\linewidth}
    \centering
    \includegraphics[keepaspectratio, scale=0.35]{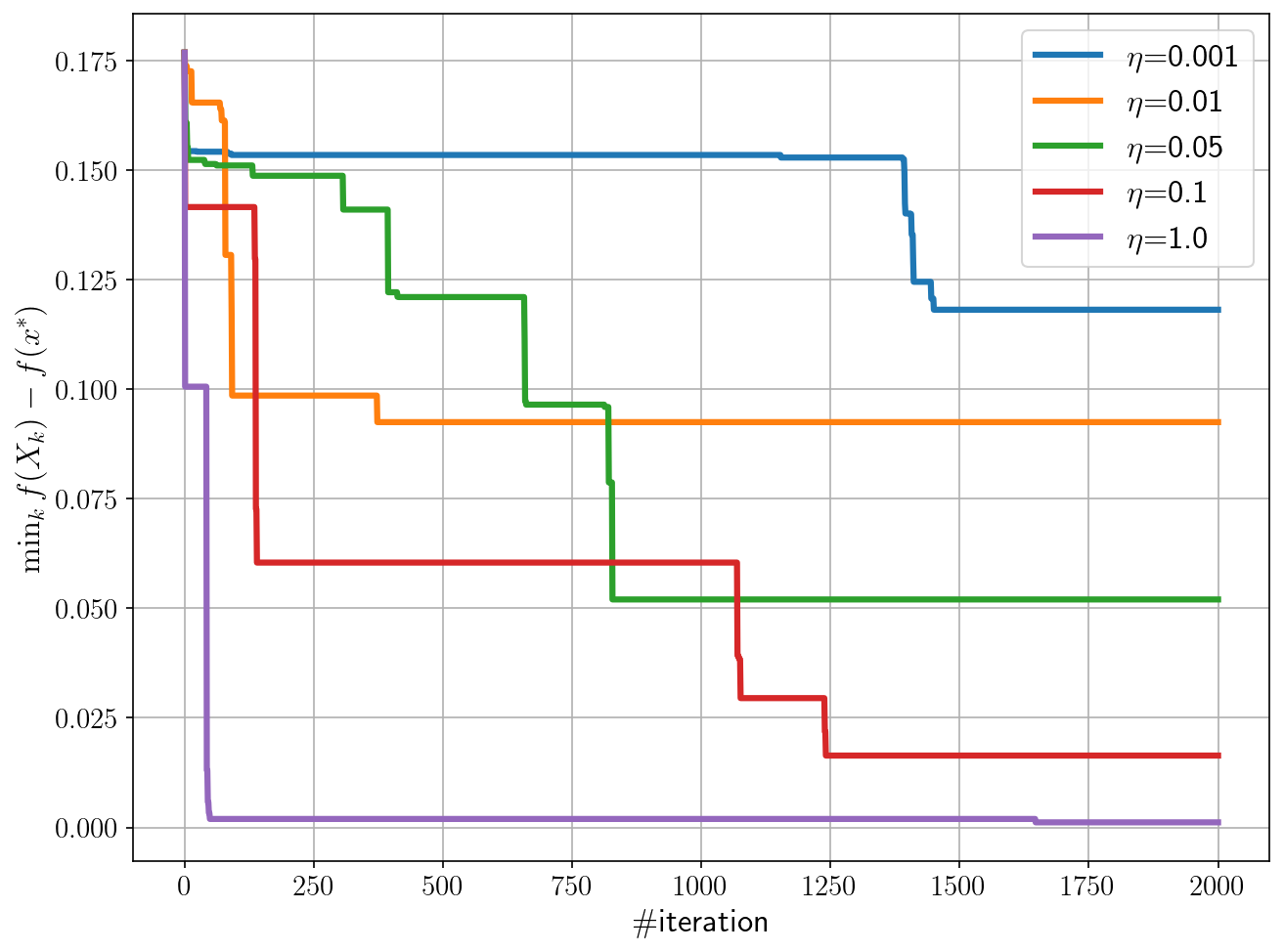}
    \subcaption{Five values of $\eta$ with $\beta=1.0$ fixed}
     \label{fig:gm-etas}
    \end{minipage} 
        \begin{minipage}[b]{0.45\linewidth}
    \centering
    \includegraphics[keepaspectratio, scale=0.32]{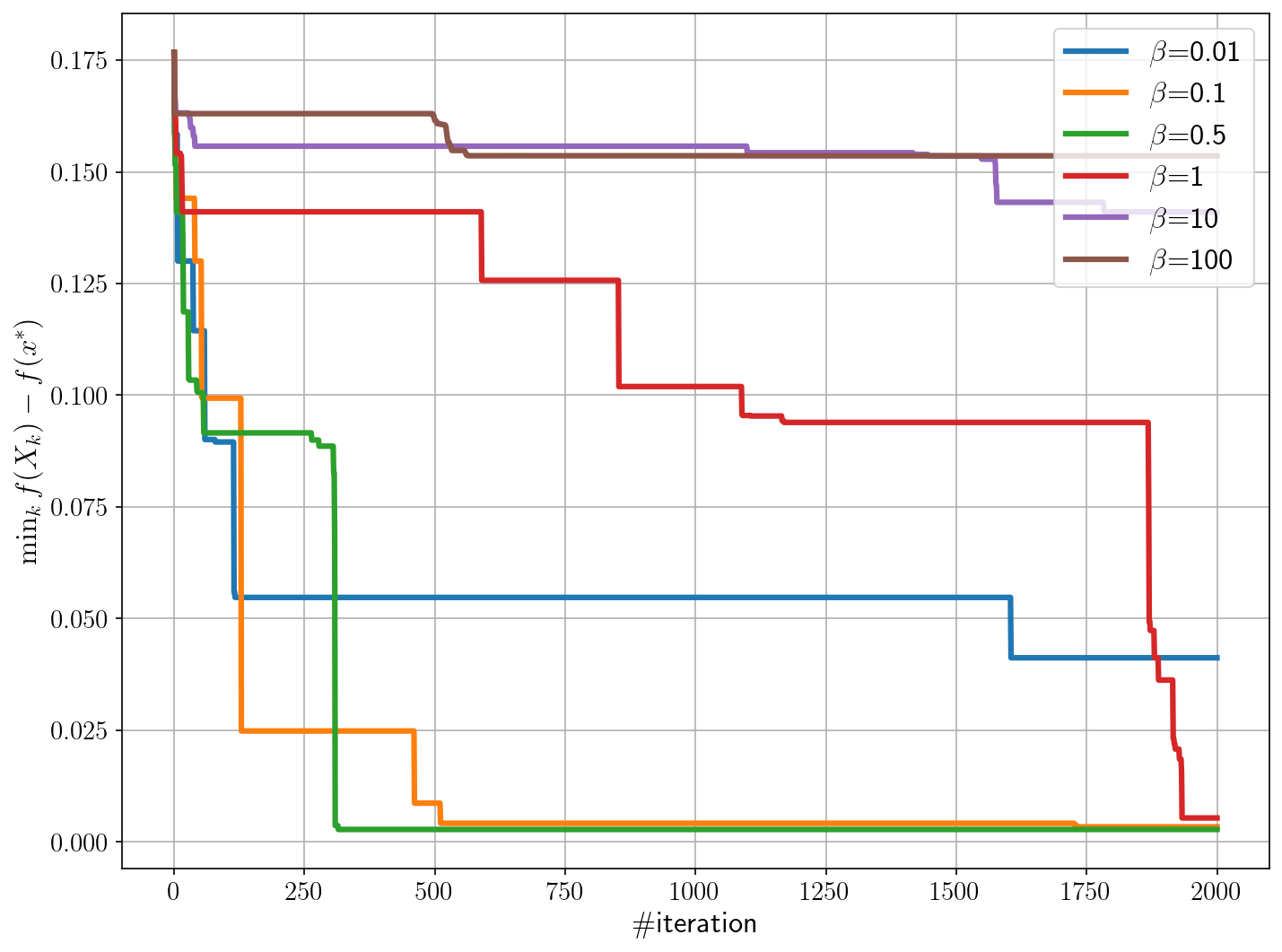}
    \subcaption{Five values of $\beta$ with $\eta=0.05$ fixed}
     \label{fig:gm-betas}
    \end{minipage} \\ ~ \\
    \begin{minipage}[b]{0.45\linewidth}
    \centering
    \includegraphics[keepaspectratio, scale=0.4]{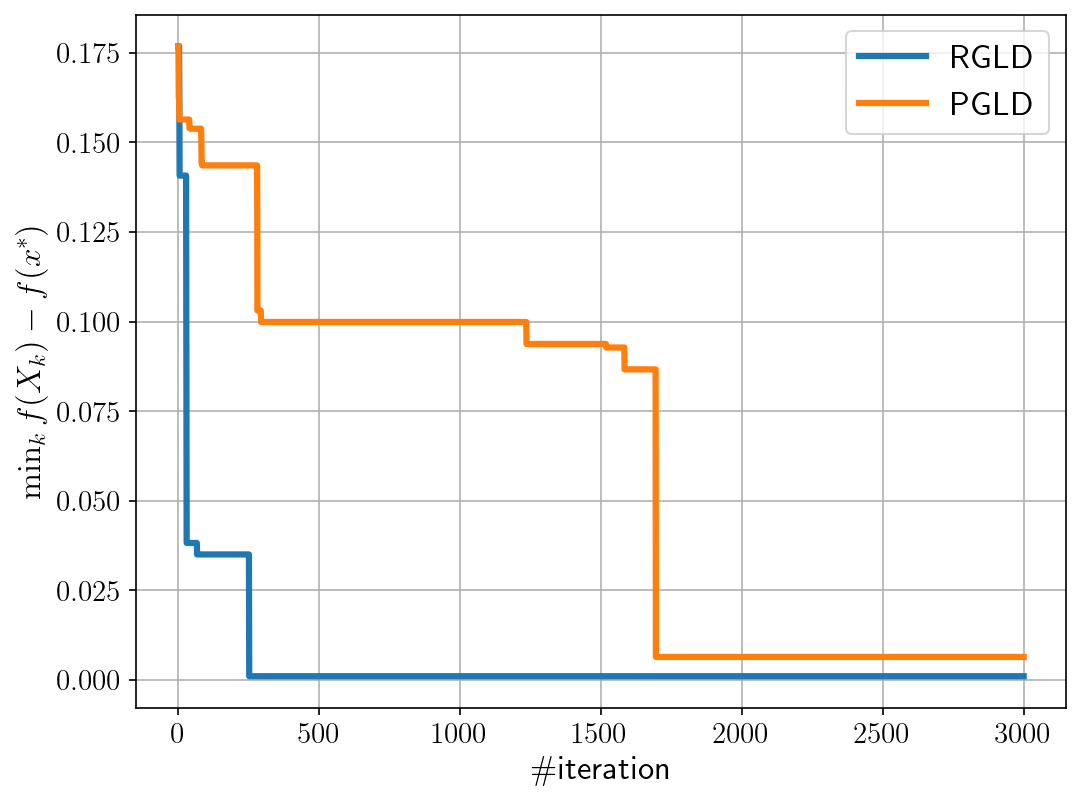}
    \subcaption{Comparison between RGLD and PGLD with $\eta=0.05$ and $\beta=1.0$ fixed}
     \label{fig:reflection-vs-projection}
      \end{minipage}
%
%
  \caption{(a)(b) Typical performance of RGLD for different values of the parameters $\eta$ and $\beta$; (c) A comparison between RGLD and PGLD in the two-dimensional problem.}
\label{fig:twodim-sensitivity}
\end{figure}

\subsection{Low-rank matrix factorization}\label{subsection low-rank matrix factorization}

Next, we examine a low-rank matrix factorization problem that has received considerable attention in the context of machine learning.
Given a matrix $X$ in $\mathbb{R}^{m\times n}$, one tries to find a pair of matrices $U\in\mathbb{R}^{m\times \ell}$ and $V\in\mathbb{R}^{n\times \ell}$ such that $X\approx UV^T$, where the column dimension $\ell$ is assumed to be much lower than the row dimensions $m$ and $n$.
If both factors $U$ and $V$ could be set arbitrarily, there would be infinitely many choices of the matrices $U$ and $V$ that yield the same matrix product.
Hence, in this experiment, we examine the following constrained problem: 
\begin{align*}
    \min_{U,V}&\quad \|X-UV^T\|_F\\
    \text{s.t.}&\quad u_1\leq \|U^i\|_2 \leq u_2,\quad i\in \{1,\cdots,m\},\\
               &\quad v_1\leq \|V^j\|_2\leq v_2,\quad j\in \{1,\cdots,n\},
\end{align*}
where $\|\cdot\|_F$ denotes the Frobenius norm, and $u_1$, $u_2$, $v_1$ and $v_2$ are suitable positive constants with $u_1<u_2$ and $v_1<v_2$.
Such $l_2$-norm constraints for matrix recovery have been discussed in \cite{haeffele2019structured,yang2019nonnegative}, backed up with promising potential for feature extraction and data representation in various fields of application, such as face recognition \cite{xu2016new}, hyperspectral unmixing \cite{zhang2018bilateral} and biomedical signal processing \cite{wang2016nmf}.

\begin{figure}[tb]
\begin{minipage}{0.38\hsize}
    \centering
    \includegraphics[keepaspectratio, width=65mm]{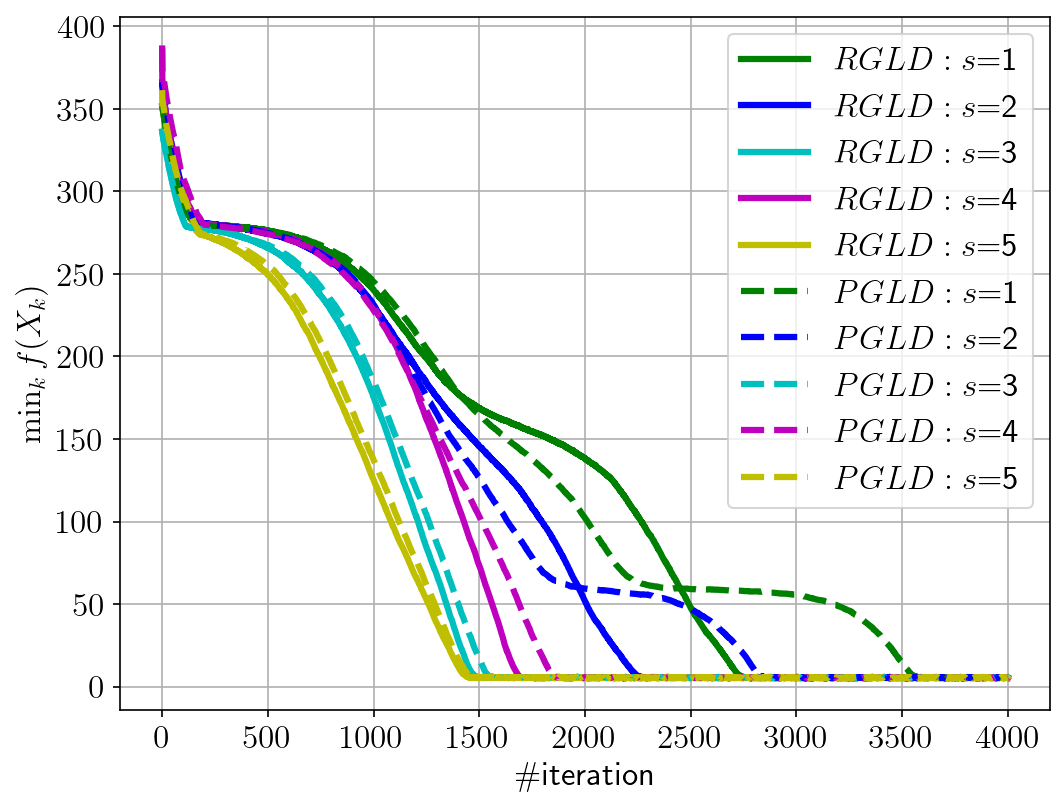}
    \subcaption{dim$=100$ and $\ell=3$.}
    \label{fig:lowrand100_3}
\end{minipage}
\hspace{1.5cm}
\begin{minipage}{0.38\hsize}
   \centering
    \includegraphics[keepaspectratio, width=65mm]{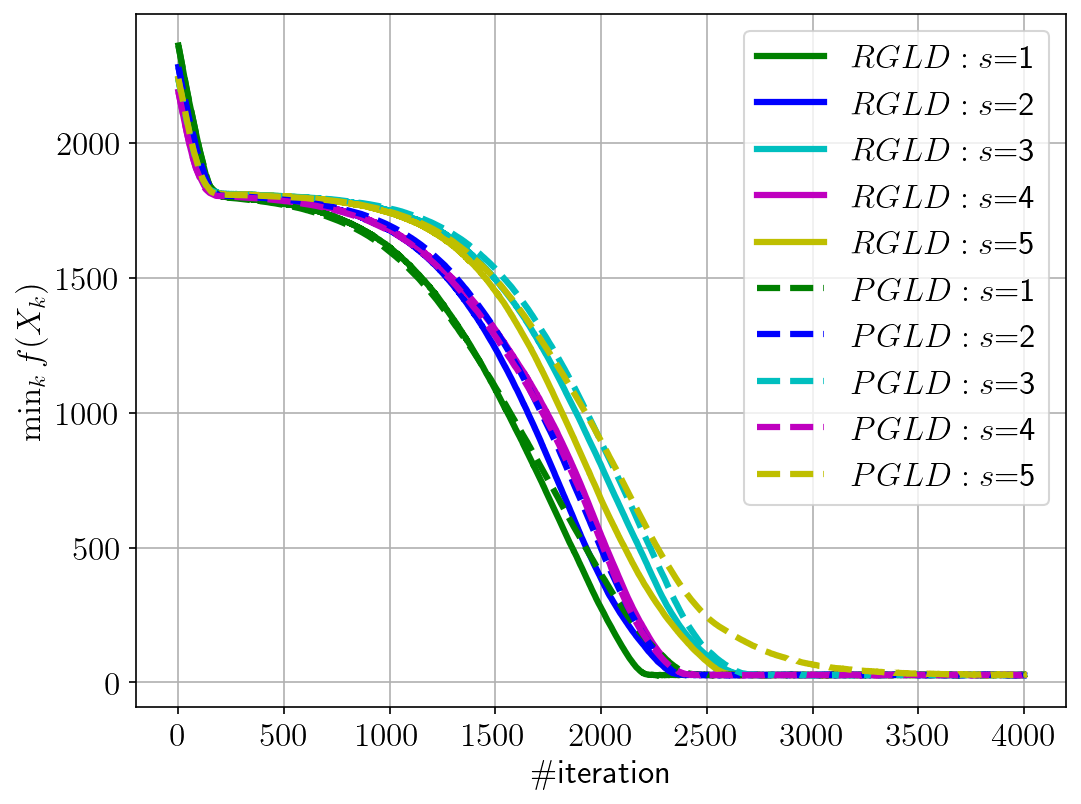}
    \subcaption{dim$=500$ and $\ell=5$.}
    \label{fig:lowrand500_5}
\end{minipage} 
\\ ~ \\ ~\\
\centering
\begin{minipage}{0.38\hsize}
  \centering
   \includegraphics[keepaspectratio, width=65mm]{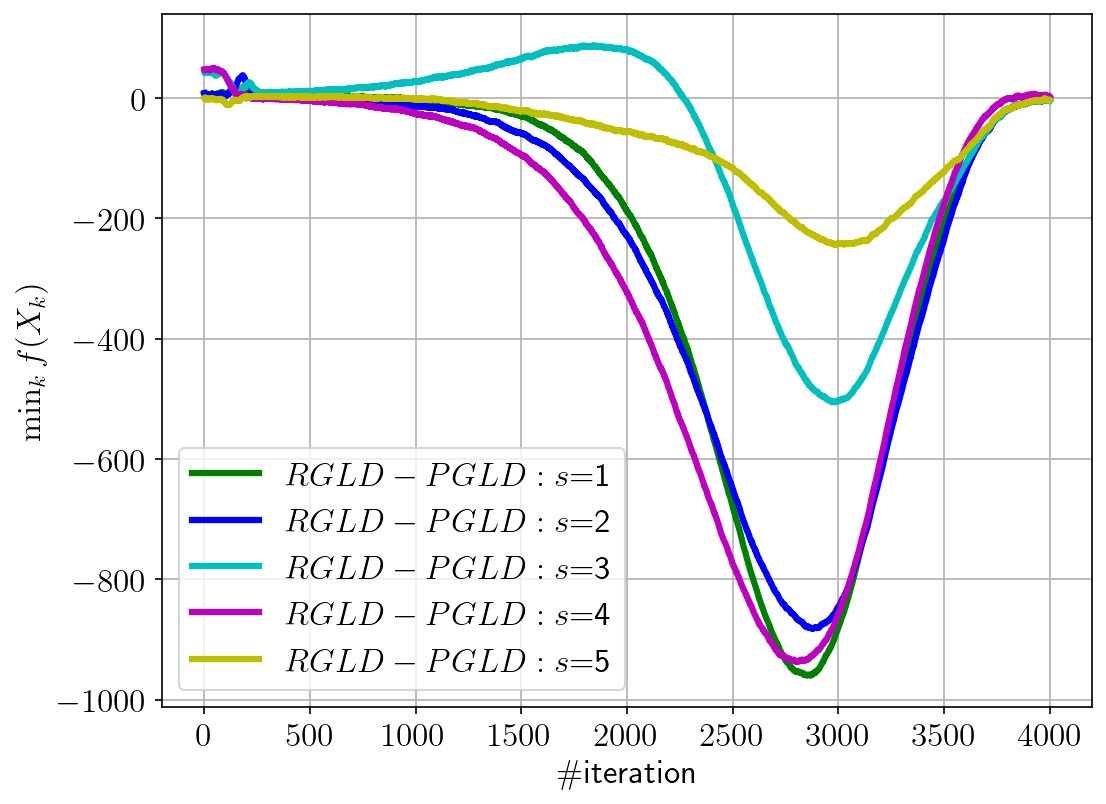}
   \subcaption{dim$=1000$ and $\ell=10$.}
   \label{fig:lowrand1000_10}
\end{minipage}
\caption{Typical iterates of RGLD and PGLD for $l_2$ regularized low-rank matrix factorization with $\eta=10^{-3}$ and $\beta=10^2$.
The initial states are randomly generated with seed $s$, as indicated.}
\label{fig:lmf-error}
\end{figure}

We examine RGLD and PGLD for three parameter settings
(a) $m=n=100$ and $\ell=3$,
(b) $m=n=500$ and $\ell=5$, and
(c) $m=n=1000$ and $\ell=10$, with five different initial states corresponding to the five different random seeds.
In every experiment, we set $u_1=v_1=0.9\sqrt{m}$ and $u_2=v_2=1.1\sqrt{m}$.
Now, Figures \ref{fig:lowrand100_3} and \ref{fig:lowrand500_5} present the function values along the run of RGLD and PGLD with initial states in common.
For better comparison, we plot gaps between those function values in Figure \ref{fig:lowrand1000_10}, the function value of RGLD minus that of PGLD, to be exact.
That is, a negative value indicates that RGLD achieves a smaller objective value than PGLD in the iteration.
In short, RGLD converges faster than PGLD with every initial state that we have examined here and proves effective enough in practice.

\section{Conclusion}
In the present work, we have derived a sub-linear convergence rate of RGLD in solving smoothly constrained problems.
In addition, the convergence rate, faster than that of \cite{Lamperski}, has been obtained by employing the reflection and making use of a direct estimate of the discretization error on the basis of the Poisson equation with the Neumann boundary condition.
The order of the spectral gap remains unclear as the smallest eigenvalue problems with the Neumann boundary condition have not been studied much.
The smoothness of the boundary, which is essential in the present analysis, is not always satisfied in many applications, whereas a non-smooth domain can be made so by a suitable smoothing operation.

We close this work by highlighting some future research directions.
It is certainly ideal that the thickness of the domain can be quantified in one way or another so as to set the step size in a more objective manner.
The scope of the algorithm, in its current form, may not sufficiently accommodate desired breadth due to the smoothness conditions outlined in Assumption \ref{assumption}.
It certainly holds significant importance in loosening these technical constraints to encompass more general unbounded and non-smooth constraints. 
One promising avenue is to introduce a dissipativity condition \cite{Raginsky, Xu} with the aid of a suitable theoretical background \cite{barbu2005neumann}, in place of the boundedness assumption of the domain.
It would also be of great interest to extend our results to the stochastic gradient oracle version.
All those would be interesting future directions of research deserving of their own separate investigation.


\vspace{1em}
\noindent {\bf Data availability statement}

\noindent All data in the paper are available from the corresponding author upon a reasonable request.
There is no conflict of interest in writing the paper.

\begin{appendices}

\section{Related studies on GLD} \label{sec:related-work}
GLD, in a broad sense including its stochastic variants (SGLD), has long been studied, dating back to \cite{geman1986diffusions}, from a wide variety of perspectives with a view towards many fields of application ever since, such as Bayesian learning \cite{Welling} to name just a few.
In particular, the convergence of GLD has been investigated quite intensively, for instance, that to the global minimum of the objective functions (\cite{Chiang,Gelfand} among many others), whereas to local optima \cite{Zhang}.

Built upon the established probabilistic analyses of GLD and SGLD,  convergence rates to global optima are derived in \cite{Raginsky,vempala2019rapid,Xu} 
for unconstrained non-convex problems in the non-asymptotic regime.
To our knowledge, the best convergence rate up to date is can be found in \cite{zou2021faster}.
SGLD is combined with stochastic variance-reduced gradient (SVRG)  \cite{chen2021approximation}, while a non-smooth objective function is investigated in \cite{durmus2018efficient} using the Moreau envelope.

Among many studies in this line of research, the work \cite{Bubeck} can be thought of as the first development on projected gradient Langevin dynamics (PGLD) for constrained sampling.
It aims at sampling from a log-concave distribution, which translates into convex optimization in our context.
Sampling from a constrained log-concave distribution is also investigated in \cite{brosse2017sampling} based on a proximal operation.
Such proximal-type algorithms have also been constructed in \cite{pereyra2016proximal,salim2019stochastic,salim2020primal}.
Moreover, mirrored Langevin dynamics is proposed in \cite{hsieh2018mirrored,zhang2020wasserstein} as a way to handle constraints, in which the mirror map is employed to transform the target distribution into an unconstrained one.
Despite this approach is valid in many existing methods, it is not a trivial task to find an appropriate mirror map in practice.
The developments \cite{patterson2013stochastic,wang2020fast} can also be considered to be related studies in the realm of constrained optimization when the feasible region forms a Riemannian manifold, such as the Stiefel or Grassmann manifold.

A GLD algorithm is proposed in \cite{parpas2006linearly} for the first time with convergence to the global optimum for linearly constrained non-convex problems.
Recently, a non-asymptotic convergence rate of PGLD is derived in \cite{Lamperski} for convex constrained non-convex problems, where reflected stochastic differential equations (RSDE) play a central role in its discretization and theoretical developments, just as in our analysis.

In probability theory, the existence and uniqueness of the solution of the reflected diffusion is studied \cite{lions1984stochastic,Tanaka} in the name of the Skorokhod problem.
The regularity of the solution and the existence of a unique stationary distribution of RSDEs are investigated, respectively, in \cite{andres2011pathwise,dupuis1999convex} and \cite{banerjee2020parameter,freidlin2016functional}.
RSDEs have also been investigated from the perspective of partial differential equations, for instance, in \cite{barbu2005neumann,cerrai1998some,miranda1970}.
The discretization error in path generation of RSDEs has been a key component in connecting the analysis on SDEs with the convergence rate of optimization algorithms \cite{cattiaux2017invariant,ding2008numerical,Leimkuhler}.


\section{Proof of Lemma \ref{lem:unique-invariant-measure}}
\label{sec:uniqueness}

In a similar manner to unconstrained problems, the existence and uniqueness can be derived for the stationary distribution of the stochastic differential equation \eqref{eq:rsde} (see, for instance, \cite{Bubeck}).
For the sake of completeness, we here outline its derivation.

Thanks to Assumption \ref{assumption}, there exists a stationary distribution of the continuous-time limit $\{X(t):\,t\ge 0\}$ uniquely.
Thus, it suffices to show that $\pi$ is stationary, which holds true if and only if $\int_K \mathcal{A}u d\pi = 0$
for all $u\in \{h \mid \langle \grad h(x),v(x)\rangle = 0,\,x\in\partial K\}$, where $\mathcal{A}$ is the infinitesimal generator defined in Section \ref{subsection main result}.
Then, one can derive
\begin{align*}
    \int_K \mathcal{A}u d\pi &= \int_K \left(\frac{1}{\beta}\Delta u(x) - \langle \grad f(x), \grad u(x)\rangle\right)\frac{\exp(-\beta f(x))}{\int_K\exp(-\beta f(y))dy}dx\\
    &= \int_K \div(\frac{\exp(-\beta f(x))\grad u(x)}{\beta\int_K\exp(-\beta f(y))dy})dx\\
    &= \int_{\partial K} 
    \frac{\exp(-\beta f(x))}{\beta\int_K\exp(-\beta f(y))dy}\langle \grad u(x), v(x)\rangle dx = 0,
\end{align*}
due to \textcolor{black}{the definitions of $\pi$ and $\mathcal{A}$} and the divergence theorem.

\section{Proof of Lemma \ref{lem:discretization-error}}
\label{sec:detailed-proof}

We take a few steps for deriving this key result.



\subsection{Error Analysis via Taylor Expansion}

Here, we aim to establish a relationship between the Taylor expansion of the solution \eqref{eq:neumann-problem} and the error analysis of the stationary distribution approximation by following the discussion of \cite{Leimkuhler}.
First, 
it is known \cite[Section 4.4]{Leimkuhler} that the solution $u$ on $\overline{K}$ can be extended to a function on $\overline{K\cup K_{-r}}$ with the degree of smoothness unchanged under Assumption \ref{assumption}-\ref{asmp:boundary-smoothness}.
Hereafter, we let $\widetilde u$ denote such a function that extends the solution $u$ to $\overline{K\cup K_{-r}}$, that is, $\widetilde{u}$ is as smooth as $u$ on $\overline{K\cup K_{-r}}$ with $\widetilde u(x)=u(x)$ for all $x \in \overline K$.
We note that $\widetilde{u}$ and its derivatives up to fourth order are uniformly bounded on the interior $K\cup K_{-r}$.
We denote by $D^i u$ the $i$-th derivative of $u$ with a linear operator on $\overline K$ based on the Riesz representation theorem.

We discuss the Taylor expansion of $u(X_{k+1})$ of \eqref{eq:neumann-problem} around the previous solution $X_k$.
For ease on notation, we write, in accordance with \cite{Leimkuhler}, $u_k:=u(X_k)$, $u'_k:=\widetilde u(X'_k)$, $v_k^{\pr} :=v(X_k^{\pr})$, $f_k:=f(X_k)$, and $X_k^{\pr}:=\mathcal{P}_K(X'_k)$, and continue to denote by $v(x)$ the outer unit normal vector at $x$ in the boundary $\partial K$, as in the main text.
Recall that the states $X_k$ and $X'_k$ reside in and strictly outside the feasible region $\overline{K}$ (Algorithm \ref{def:rgld}), resulting in \textcolor{black}{$\widetilde{u}(X_k)=u(X_k)$}, while $X_k'$ cannot be evaluated with $u$ but needs to be with the extension $\widetilde{u}$.
Following this notation, we reroute the increment as
    $u_{k+1} - u_k = (u_{k+1} - u_{k+1}') +(u_{k+1}'-u_k)$.
As for the first term, since $X_{k+1}^{\pr}$ is the midpoint between $X'_{k+1}$ and $X_{k+1}$ in accordance with \eqref{def:projection}, it holds that
\[
X_{k+1} = \mathcal{R}_K(X'_{k+1}) = X'_{k+1} - 2 r_{k+1} v_{k+1}^{\pr} = X^{\pr}_{k+1} - r_{k+1} v_{k+1}^{\pr},
\]
and
\[
    X'_{k+1} = X'_{k+1} - r_{k+1} v_{k+1}^{\pr} + r_{k+1} v_{k+1}^{\pr} = X^{\pr}_{k+1} + r_{k+1} v_{k+1}^{\pr},
\]
where $r_k:=\|X_k^{\pr}-X_k\|$ for $k\in \mathbb{N}$.
Hence, by performing the third-order Taylor expansion and applying the mean value theorem to the third derivative term, it holds that for suitable constants $\alpha_1$ and $\alpha_2$ in $[0,1]$,
\begin{align}
        u_{k+1}-u_{k+1}' &=u(X_{k+1}^{\pr}-r_{k+1} v_{k+1}^{\pr}) - \widetilde u(X_{k+1}^{\pr}+r_{k+1} v_{k+1}^{\pr})\nonumber\\
        &= \left(u_{k+1}^{\pr}-r_{k+1}Du_{k+1}^{\pr}\left(v_{k+1}^{\pr}\right)+\frac{r_{k+1}^2}{2}D^2u_{k+1}^{\pr}\left(v_{k+1}^{\pr},v_{k+1}^{\pr}\right)\right.\nonumber\\
        &\qquad \left.-\frac{r_{k+1}^3}{6}D^3 \widetilde u(X_{k+1}^{\pr}-\alpha_1r_{k+1} v_{k+1}^{\pr})\left(v_{k+1}^{\pr},v_{k+1}^{\pr},v_{k+1}^{\pr}\right)\right)\nonumber\\
        &\qquad - \left(u_{k+1}^{\pr}+r_{k+1}Du_{k+1}^{\pr}\left(v_{k+1}^{\pr}\right)+\frac{r_{k+1}^2}{2}D^2u_{k+1}^{\pr}\left(v_{k+1}^{\pr},v_{k+1}^{\pr}\right)\right.\nonumber\\
        &\qquad \left.+\frac{r_{k+1}^3}{6}D^3 \widetilde u(X_{k+1}^{\pr}+\alpha_2r_{k+1} v_{k+1}^{\pr})\left(v_{k+1}^{\pr},v_{k+1}^{\pr},v_{k+1}^{\pr}\right)\right)\nonumber\\
        &= -\frac{r_{k+1}^3}{6}D^3 \widetilde u(X_{k+1}^{\pr}-\alpha_1r_{k+1} v_{k+1}^{\pr})\left(v_{k+1}^{\pr},v_{k+1}^{\pr},v_{k+1}^{\pr}\right)\nonumber\\
        & \qquad- \frac{r_{k+1}^3}{6}D^3 \widetilde u(X_{k+1}^{\pr}+\alpha_2r_{k+1} v_{k+1}^{\pr})\left(v_{k+1}^{\pr},v_{k+1}^{\pr},v_{k+1}^{\pr}\right)\label{eq:first-taylor}\\
        &\ge -\frac{r_{k+1}^3}{3}\norm{D^3 \widetilde u},\nonumber
\end{align}
where we have applied $Du_{k+1}^{\pr}(v_{k+1}^{\pr})=\langle \grad u(X^{\pr}_{k+1}), v_{k+1}^{\pr}\rangle=0$, due to \eqref{eq:neumann-problem}.
The last inequality holds true by the definition of the operator norm and since $\|v_{k+1}^{\pr}\|$ only takes values in $\{0,1\}$.
In a similar manner, with a constant $\alpha_3$ in $[0,1]$, we obtain
\begin{align*}
        u'_{k+1}-u_k
        & = \widetilde u\left(X_k-\eta\grad f_k+\sqrt{\frac{2\eta}{\beta}}\xi_{k+1}\right) - u_k\\
        &= Du_k\left(-\eta\grad f_k+\sqrt{\frac{2\eta}{\beta}}\xi_{k+1}\right)\\
        &\quad + \frac{1}{2}D^2u_k\left(-\eta\grad f_k+\sqrt{\frac{2\eta}{\beta}}\xi_{k+1},-\eta\grad f_k+\sqrt{\frac{2\eta}{\beta}}\xi_{k+1}\right)\\
        &\quad + \frac{1}{6}D^3\widetilde u\left(X_k + \alpha_3\left(-\eta\grad f_k+\sqrt{\frac{2\eta}{\beta}}\xi_{k+1}\right)\right)\left(-\eta\grad f_k+\sqrt{\frac{2\eta}{\beta}}\xi_{k+1},\right.\\
        &\quad \left.-\eta\grad f_k+\sqrt{\frac{2\eta}{\beta}}\xi_{k+1},-\eta\grad f_k+\sqrt{\frac{2\eta}{\beta}}\xi_{k+1}\right),
\end{align*}
and 
\begin{align}
        \mathbb{E}\left[u'_{k+1}-u_k\right] 
        & = \mathbb{E}\left[\eta \left(f_k-\mathbb{E}_\pi f\right)+\frac{\eta^2}{2}D^2u_k\left(\grad f_k,\grad f_k\right)\right.\nonumber\\
        &\quad \left.+\frac{1}{6}D^3 \widetilde u\left(X_k + \alpha_3\left(-\eta\grad f_k+\sqrt{\frac{2\eta}{\beta}}\xi_{k+1}\right)\right)\left(-\eta\grad f_k+\sqrt{\frac{2\eta}{\beta}}\xi_{k+1},\right.\right.\nonumber\\
        &\quad \left.\left.-\eta\grad f_k+\sqrt{\frac{2\eta}{\beta}}\xi_{k+1},-\eta\grad f_k+\sqrt{\frac{2\eta}{\beta}}\xi_{k+1}\right)\right],\label{eq:second-taylor}
\end{align}
due to \eqref{infinitesimal generator} and $\mathbb{E}[\xi_{k+1}]=0$.
By applying \eqref{eq:first-taylor} and \eqref{eq:second-taylor} to the equation $\mathbb{E}[u_{k+1}'-u_k] = \mathbb{E}[u_{k+1}-u_k] - \mathbb{E}[u_{k+1}-u_{k+1}']$, we obtain 
\begin{align*}
    &\mathbb{E}\left[\eta \left(f_k-\mathbb{E}_\pi f\right)+\frac{\eta^2}{2}D^2u_k\left(\grad f_k,\grad f_k\right)\right.\\
    &\qquad \left.+\frac{1}{6}D^3 \widetilde u\left(X_k + \alpha_3\left(-\eta\grad f_k+\sqrt{\frac{2\eta}{\beta}}\xi_{k+1}\right)\right)\left(-\eta\grad f_k+\sqrt{\frac{2\eta}{\beta}}\xi_{k+1},\right.\right.\\
    &\qquad \left.\left.-\eta\grad f_k+\sqrt{\frac{2\eta}{\beta}}\xi_{k+1},-\eta\grad f_k+\sqrt{\frac{2\eta}{\beta}}\xi_{k+1}\right)\right]\\
    & \qquad \qquad \le \E{u_{k+1}-u_k} + \frac{r_{k+1}^3}{3}\norm{D^3 \widetilde u}.
\end{align*}
Since $\|\xi_{k+1}\|=\sqrt{d}$ and $\|\grad f\|$ is uniformly bounded over the domain $K$, we get
\begin{equation}
\label{eq:error-estimation}
    \E{\eta \left(f_k-\mathbb{E}_\pi f\right)}\leq \E{u_{k+1}-u_k}
    + \frac{r_{k+1}^3}{3}\norm{D^3 \widetilde u} + C\eta^2\norm{D^2 \widetilde u} + C\norm{D^3 \widetilde u}\left(\eta+\sqrt{d\eta/\beta}\right)^3.
    \end{equation}
To obtain the upper bounds of the right-hand side, we employ the following result, which we derive after completing the proof of Lemma \ref{lem:discretization-error} (Appendix \ref{sec:ProofOfLemmaDerivatives}) to maintain the flow of the paper.

\begin{lem}
\label{lem:derivative-bounds}
Let $i\in \{0,1,2,3\}$ and let $r$ be the supremum of positive numbers such that the projection of all $x$ in $K_{-r}$ is uniquely defined.
There exists a positive constant $C$ such that $\|D^i\widetilde u(x)\|\le C\beta^{i/2}/\lambda_*$ for all $x\in \overline{K\cup K_{-r}}$.
\end{lem}

To continue Lemma \ref{lem:discretization-error}, we apply Lemma \ref{lem:derivative-bounds} to the inequality \eqref{eq:error-estimation} to obtain
\[
    \mathbb{E}\left[\eta\left(f(X_k)-\mathbb{E}_\pi f\right)\right] 
    \le \mathbb{E}\left[u(X_{k+1})-u(X_k)\right] + \frac{C}{\lambda_*}\left(\beta^{3/2}r_{k+1}^3+\beta\eta^2+\beta^{3/2}\eta^3+(d\eta)^{3/2}\right),
\]
which, with $r_{k+1}=\mathcal{O}(\eta+\sqrt{d\eta/\beta})$ shown in \eqref{eq:proj_dist}, further yields
\[
    \mathbb{E}\left[\eta\left(f(X_k)-\mathbb{E}_\pi f\right)\right] < \mathbb{E}\left[u(X_{k+1})-u(X_k)\right] + \frac{C\eta^{3/2}\left(\beta\sqrt{\eta}+(\beta\eta+d\right)^{3/2})}{\lambda_*}.
\]
Thus, again thanks to Lemma \ref{lem:derivative-bounds}, we get
\begin{align*}
    \frac{1}{N}\sum_{k\in [N]} \mathbb{E}\left[f(X_k)-\mathbb{E}_\pi f\right] &\le  \mathbb{E}\left[\frac{u(X_N)-u(X_0)}{N\eta}\right]+\frac{C\sqrt{\eta}\left(\beta\sqrt{\eta}+(\beta\eta+d)^{3/2}\right)}{\lambda_*}\\
    &\le \frac{C}{\lambda_*N\eta}+\frac{C\sqrt{\eta}\left(\beta\sqrt{\eta}+(\beta\eta+d)^{3/2}\right)}{\lambda_*},
\end{align*}
which derives Lemma \ref{lem:discretization-error}.

\subsection{Proof of Lemma \ref{lem:derivative-bounds}}
\label{sec:ProofOfLemmaDerivatives}

Finally, we give the proof of Lemma \ref{lem:derivative-bounds}. 
Towards this end, we show the following Lemma~\ref{lem:exponential-convergence-derivatives}
on the function $\phi(t, x):=\mathbb{E}\left[f(X^x(t))-\mathbb{E}_\pi f\right]$.

\begin{lem}
    \label{lem:exponential-convergence-derivatives}
Let $i\in \{1,2,3\}$.
There exists a positive constant $C$ such that 
    $\|D^i\phi(t,x)\|\le C\beta^{i/2}e^{-\lambda_* t}$ for all $t\geq 0$ and $x\in \overline K$.
\end{lem}

Integrating the inequality in $t$ together with $u(x)=-\int_0^\infty \phi(t, x) dt$ for $x\in \overline K$ and
extending the range of $x$ from $\overline K$ to $\overline{K\cup K_{-r}}$ as discussed in the beginning of this section
lead to the claim of Lemma~\ref{lem:derivative-bounds}.
Therefore, it now suffices to show Lemma \ref{lem:exponential-convergence-derivatives}.
{To this end, we employ the following two estimates, due to \cite[Section 3.2]{refId0} and \cite[Lemma 3.7.3]{freidlin2016functional}, respectively.
We note that the second one does not stand alone by itself for our present purpose but helps the first one in combination.}

\begin{lem} 
    \label{lem:estimates}
    (i) 
    Let $i\in \{1,2,3\}$, 
    $g\in \mathcal{C}^5(\overline{K})$,  
    and define $\Phi(t,x):=\mathbb{E}[g(X^x(t))]$.
    There exists a positive constant $C$ such that $\|D^i\Phi(t,x)\|\le C (\beta/t)^{i/2}\|g\|_\infty$
for all $t\in (0,1]$ and $x\in \overline K$.

\noindent (ii) There exists a positive constant $C$ such that $|\phi(t,x)| \le Ce^{-\lambda_*t}$ for all $t\geq 0$ and $x\in \overline K$.
\end{lem}




To continue, fix an initial state $x$ in the domain $K$ and define the (random) set $S:=\{t\geq 0:\,X^x(t)\in\partial K\}$ and the function $r(t):=\sup(S\cap[0, t])$.
For $h \in \mathbb{R}^d$ and $t \geq 0$, let $\{Y_t^h:\,t\ge 0\}$ be the solution to
\[
Y^h_t=
\begin{cases}
    h - \int_0^t \Delta f(X^x(s)) Y^h_s ds,& \text{if}\: t<\inf S,\\
    \mathcal{P}_{X^x(r(t))}(Y^h_{r(t)^-})-\int_{r(t)}^t \Delta f(X^x(s)) Y^h_sds,&\text{if}\: t\geq \inf S,
\end{cases}
\]
where $Y^h_{s^-}$ denotes the left limit $\lim_{t \uparrow s^-}Y^h_t$ and $\mathcal{P}_x(z):= z - \langle z, v(x)\rangle v(x)$, for $z\in\mathbb{R}^d.$
For ease of notation, we have subscripted the (continuous) time index of the solution $\{Y_t^h:\,t\ge 0\}$ (as opposed to $\{X^x(t):\,t\ge 0\}$). 
It is known (for instance, \cite{andres2011pathwise}) that $Y_t^h$ is the first derivative (in probability) of $X^x(t)$ with respect to $x$ along the direction $h\in\mathbb{R}^d$.
Next, we denote by $\{Z_t^{h,k}:\,t\ge 0\}$ the solution to 
\[
Z_t^{h,k} =
\begin{cases}
    -\int_0^t \Delta f(X^x(s)) Z_s^{h,k}+D^3 f(X^x(s))(Y_s^h,Y_s^k)ds,& \text{if } t<\inf S,\\
    -\int_{r(t)}^t \Delta f(X^x(s)) Z_s^{h,k}+D^3 f(X^x(s))(Y_s^h,Y_s^k)ds,& \text{if } t\geq\inf S.
\end{cases}
\]
In a similar manner to the first derivative, it is known (for instance, \cite{cerrai1998some}) that $Z_t^{h,k}$ is the second derivative (in probability) of $X^x(t)$ with respect to $x$ along the directions $h$ and $k$.
Finally, let $\{V_t^{l,m,n}:\,t\ge 0\}$ denote the solution to
\[
V_t^{l,m,n} =
\begin{cases}
     -\int_0^t (\Delta f(X^x(s))V_s^{l,m,n} + D^4 f(X^x(s))(Y_s^{l},Y_s^{m},Y_s^{n})\\
     \qquad\qquad+ \frac{1}{2}\sum_{\sigma\in \mathcal{S}_3}D^3 f(X^x(s))(Y_s^{\sigma(l)}, Z_s^{\sigma(m),\sigma(n)})) ds,& \text{if } t<\inf S,\\
    -\int_{r(t)}^t (\Delta f(X^x(s))V_s^{l,m,n} + D^4 f(X^x(s))(Y_s^{l},Y_s^{m},Y_s^{n})\\ 
    \qquad\qquad+ \frac{1}{2}\sum_{\sigma\in \mathcal{S}_3}D^3 f(X^x(s))(Y_s^{\sigma(l)}, Z_s^{\sigma(m),\sigma(n)}) )ds,& \text{if } t\geq\inf S,
\end{cases}
\]
where $\mathcal{S}_3$ denotes the set of all permutations of three elements.
Just as for the first and second derivatives, $V_t^{l,m,n}(x)$ is the third derivative (in probability) of $X^x(t)$ with respect to $x$ along the directions $l$, $m$ and $n$.

We are now ready to derive Lemma \ref{lem:exponential-convergence-derivatives} with the aid of
Lemma \ref{lem:estimates}.
  If $t\geq 1$, then
  applying Lemma \ref{lem:estimates} (i) with $t \leftarrow 1$ and $g(x)\leftarrow \phi(t-1, x)$ in its statement and evaluating $\norm{g}_\infty$ by Lemma \ref{lem:estimates} (ii) yields
    \[
        \norm{D^i\E{\phi(t-1,X^x(1))}} \le C\beta^{i/2}\mathrm{e}^{-\lambda_*(t-1)},
    \]
which is independent of the initial state $x$ and is thus equivalent to
\[
\norm{D^i\phi(t, x)} \le C\beta^{i/2}\mathrm{e}^{-\lambda_*(t-1)}\le C\beta^{i/2}\mathrm{e}^{-\lambda_*t},
\]
which proves the claim. 

As for the situation $t < 1$, we start by showing that $Y_t$, $Z_t$ and $V_t$ are almost surely bounded for all $t\in (0,1]$.
For $t< \inf S$, we have
    \begin{align*}
        \norm{Y_t^h} &\leq \norm{h} + \int_0^t \norm{\Delta f(X^x(s)) Y_s^h} ds
        \leq \norm{h} + C \int_0^t \norm{Y_s^h} ds.
    \end{align*}
    Then, with the aid of Gronwall's inequality, we get
    \begin{align*}
        \norm{Y_t^h}\leq \norm{h}\exp(\int_0^t Cds) \leq \norm{h}\exp(\int_0^1 Cds) \leq \norm{h}\mathrm{e}^C \le C\norm{h}.
    \end{align*}
    If $t\ge \inf S$, it holds by recursion that
    \begin{align*}
        \norm{Y_t^h} \le C\norm{\mathcal{P}_{X^x({r(t)})}(Y_{r(t)-}^h)} \le C\norm{Y_{r(t)-}^h} \le C\norm{h}.
    \end{align*}
    As for $Z_t$, we obtain from the results for $Y_t^h$ above that
    \begin{align*}
        \norm{Z_t^{h,k}} \leq \int_0^t C\norm{Z_s^{h,k}} + C\norm{h}\norm{k} ds \leq C\norm{h}\norm{k} + \int_0^t C\norm{Z_s^{h,k}} ds.
    \end{align*}
    Again with the aid of Gronwall's inequality, we get $\|Z_t^{h,k}\|\leq C\norm{h}\norm{k}$.
    In a similar manner, it follows that $V_t^{l,m,n}\leq C\norm{l}\norm{m}\norm{n}$.
    Finally, for $i\in \{1,2,3\}$, the quantity $D^i\phi(t,x)$ can be represented, as follows:
    \begin{align*}
        D\phi(t,x).h &= \E{Df(X^x(t)).Y_t^h(x)},\\
        D^2\phi(t,x)\left(h, k\right) &= \E{D^2f(X^x(t))\left(Y_t^h, Y_t^k\right)+Df(X^x(t)).Z_t^{h,k}},\\
        D^3\phi(t,x)\left(l,m,n\right)&= \E{D^3f(X^x(t))\left(Y_t^l,Y_t^m,Y_t^n\right)+Df(X^x(t)).V_t^{l,m,n}}\\
        &\quad +\E{\frac{1}{2}\sum_{\sigma\in\mathcal{S}_3}D^2f(X^x(t))\left(Y_t^{\sigma(l)},Z_t^{\sigma(m),\sigma(n)}\right)}.
    \end{align*}
    Combining those with the above statement yields upper bounds: $\|D\phi(t,x).h\| \leq C\norm{h}$, $\|D^2\phi(t,x)\left(h, k\right)\| \leq C\norm{h}\norm{k},$ and $\|D^3\phi(t,x)\left(l,m,n\right)\|\leq  C\norm{l}\norm{m}\norm{n}$, that is, 
    the quantities $\|D^i\phi(t,x)\|$ are bounded by suitable constants for all $t\in (0,1].$
    This yields Lemma \ref{lem:exponential-convergence-derivatives}.

\end{appendices}

\bibliography{ref}

\end{document}